\documentclass[microtype]{gtpart}
\usepackage{pinlabel}

\usepackage{graphicx}
\usepackage{multirow}
\usepackage[english]{babel}
\usepackage[alphabetic]{amsrefs}
\graphicspath{ {images/} }
\usepackage[T1]{fontenc}
\usepackage{txfonts}
\usepackage{times}
\usepackage[all]{xy}
\usepackage{subfig}
\usepackage{tikz}
\usetikzlibrary{shapes,arrows,shadows}
\usetikzlibrary{decorations.markings}
\usepackage{color}
\usepackage[normalem]{ulem}
\usepackage{hyperref}
\usepackage{wrapfig}
\usetikzlibrary{arrows}
\usepackage{pb-diagram}
\usepackage{verbatim}
\usepackage{float}
\usepackage{pinlabel}

\newtheorem{thm}{Theorem}[section]

\newtheorem{prop}{Proposition}[section]
\newtheorem{corolario}{Corollary}[section]
\newtheorem{definition}{Definition}[section]

\newtheorem{remark}{Remark}[section]

\begin{document}

\title{Symmetry Type Graphs on 4-Orbit Maps}

\author[Arredondo]{John A. Arredondo}
\givenname{John A.}
\surname{Arredondo}
\email{alexander.arredondo@konradlorenz.edu.co}
\address{John A. Arredondo\\\newline Fundaci\'on Universitaria Konrad Lorenz\\\newline CP. 110231, Bogot\'a, Colombia.}
%\urladdr{http://docentes.konradlorenz.edu.co/2015/06/john-alexander-arredondo-g.html}

\author[Ram\'irez Maluendas]{Camilo Ram\'irez Maluendas}
\givenname{Camilo}
\surname{Ram\'irez Maluendas}
\email{camramirezma@unal.edu.co}
\address{Camilo Ram\'irez Maluendas\\\newline Universidad Nacional de Colombia, Sede Manizales\\\newline Manizales, Colombia.}
%\urladdr{http://www.fcen.unal.edu.co/menu/departamentos/departamento-de-matematica-y-estadistica/profesores/}

\author[Santos Guerrero]{Luz Edith Santos Guerrero}
\givenname{Luz Edith}
\surname{Guerrero Santos}
\email{luze.santosg@konradlorenz.edu.co}
\address{Luz Edith Santos Guerrero\\\newline Fundaci\'on Universitaria Konrad Lorenz\\\newline CP. 110231, Bogot\'a, Colombia.}
%\email{luze.santosg@konradlorenz.edu.co}



\begin{abstract}
 It is well known that there exist twenty two symmetry type graphs associated to $4$-orbit maps. For this ones
 we give the feasible values taken by the degree of the vertices and the number appropriate of edges in the boundary of each face of the map, by introducing the concepts of vertex type graph, face type graph and characteristic system.
\end{abstract}

\maketitle

%\begin{center}\rule{1.0\textwidth}{0.1mm} \end{center}
%\begin{abstract}
% It is well known that there exist twenty two symmetry type graphs associated to $4$-orbit maps. For this ones
% we give the feasible values taken by the degree of the vertices and the number appropriate of edges in the boundary of each face of the map, by introducing the concepts of vertex type graph, face type graph and characteristic system.
%\end{abstract}

%\textbf{Keywords:}  Pregraph, Symmetry type graph, 4-orbit maps, vertex type graph, face type graph, characteristic system.

%\textbf{AMS:}	05C30, 05C25, 52B15, 05C07. 

%\begin{center}\rule{1.0\textwidth}{0.1mm} \end{center}
%%\end{@twocolumnfalse}
%%]

%\tableofcontents

\section{Introduction}\label{sec:introduction}

The concept of map on a surface $S$, comes from the ancient idea of a map of the Earth. The surface $S$ is decomposed into countries ``faces'' where every border ``edge'' belongs to exactly two countries. The points where three or more countries are incident correspond to the vertices of the map. %In each face we can choose an interior point and define a triangulation of the surface $S$, these triangles are called the flags of the map.
For this, we  fix a point in the interior of each face, called the center of the face. Thus, in  each face, we draw a line segments with nodes the center and the vertices, which are on the boundary of each one of them, respectively. Likewise, on each edge in the boundary of each face we mark the middle point and, we draw the line segments from the center of each face to the middle points on each one its edges, respectively. This process defines a triangulation of the surface $S$, and each topological triangle is called a flag of the map.

% From these elements we can define the triplet $(v,e, f)$ called a flag of the map, which is confromed by a vertex $v$, an edge $e$ and, a face $f$, such that the edge $e$ is belonged to the border of the face $f$ and, the vertex $e$ is incident with the edge $e$.}

%(trazamos curvas desde el punto interior hasta los vertices, dichas curvas solo tienen en comun el punto interior. Luego en cada en cada arista tomamos el punto medio y nuevamente por cada punto medio en cada arista tramos una curvas que no intersactan a las ya trazadas y que tiene extremos el punto interior de la cara y la respectiva arista)

Every map has associated  its automorphism group  and a pregraph, called ``symmetry type graph'', built as  the quotient of its flag graph  under the action of the automorphism group, this object has been investigated by Cunningham \emph{at all} \cite{Maniplexes}, Kovi\v{c}  \cite{Kov},  Del Rio \cite{Fran}, and Hubard \cite{HubardTesis}. As the automorphism group of the map acts freely on the set conformed by all flags of the map, if the action defines $k$ classes, then the map is said to be a $k$-orbit map. In \cite[Chapter 8]{CoMo} H. S. M. Coxeter and W. O. J. Moser called to the $1$-orbit maps, regular maps and a class of $2$-orbit maps, irreflexible maps, also known as chiral maps. Regular and irreflexible maps have been studied widely by several authors  as \emph{e.g.}, Wilson \cite{WilSte},  D'Azevedo,  Jones and Schulte \cite{DaJonSc}, and the first and third author jointly with  Valdez  \cite{ARV} among others, because regular and chiral maps are the most symmetric ones.  Mostly, the $k$-orbit maps on surfaces are interesting for their large number of implications and because in this subject converge topics as algebraic geometric, combinatorics and topology, reason for which  it has attracted the attention of numerous researchers, see \emph{e.g.}, the work by Del Rio \cite{delRio}, Helfand \cite{HelIL}, Cunningham and Pellicer \cite{CunPell}. In this paper we focus on $4$-orbit maps. Specifically, from the symmetry type graph associated to a $4$-orbit map, we give the feasible values taken by the degree of the vertices and the number appropriate of edges in the boundary of each face of the map. 

This article is organized as follows. In Section \ref{sec:map} we introduce some elements of the theory of maps as their flags, $k$-orbit map and automorphism and monodromy groups. In Section \ref{sec:pregraph} we explore the concept of a symmetry type graph associated to a map, and we present the twenty two pregraphs, which could be symmetry type graph of any $4$-orbit map. Moreover, we describe the dual and petrial of a map.  Finally, in Section \ref{sec:charsys} we introduce the concept of characteristic system of a vertex and a face. Also, we associate suitably to each vertex and each face of any $4$-orbit map a pregraph. From these elements we summarize through a table, the feasible values taken by the degree of the vertices and the number appropriate of edges in the boundary of each face of $4$-orbit maps.

%%%%%%%%%%%%%%%%%%%%%%%%%%%%%%%%%%%%%%%%%%%%%%%%%%%%%%%%%

\section{Some review on maps}\label{sec:map}

Along this paper, the term \textbf{surface} means a connected $2$-dimensional topological real manifold with empty boundary, and it will be denoted as $S$. In particular,
the transition functions of the corresponding atlas are only required to be continuous. It is important to remark that we  require $S$ to be a compact topological space.

In this text the object \textbf{map $\mathcal{M}$ on a surface $S$} means a finite $2$-cell embedding $i:G\hookrightarrow S$  of a locally finite simple graph $G$\footnote{For us $G$ will be the geometric realization of an abstract graph.} into $S$. In other words, only a finite number of edges are incident in each vertex of $G$, the vertices  of each edge are in different vertices and the function $i$ is a topological embedding, such that each connected component of  $S\setminus i(G)$ is homeomorphic to an open disk, whose boundary is the image under $i$ of a closed finite path in $G$.

Each connected component of $S\setminus i(G)$ is called a \textbf{face} of the map $\mathcal{M}$. A \textbf{vertex $i(v)$ of the map}  is the image under $i$ of a vertex $v$ in $G$. Likewise, an \textbf{edge $i(e)$ of the map}  is the image under $i$ of  an edge $e$ in $G$. The \textbf{degree} of $i(v)$ with $v\in V(G)$ is the degree of  $v$.  The  \textbf{size} of a face $f$ of the map $\mathcal{M}$ is the number of edges conforming its boundary.

\begin{figure}[h!] 
  \centering
  \includegraphics[scale=0.25]{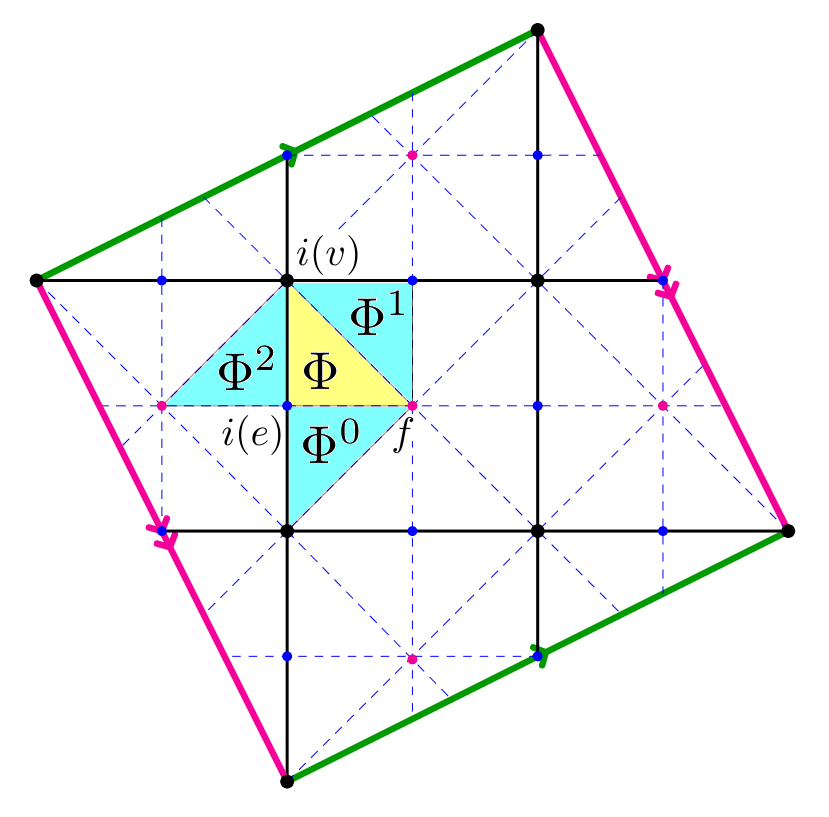}
  \caption{\emph{A map on the torus divided into flags.}}
  \label{Im:flags}
\end{figure}

A \textbf{flag} $\Phi$ of the map $\mathcal{M}$ is a triangle on the surface $S$ whose vertices are: a vertex $i(v)$, the ``midpoint'' of an edge $i(e)$ incident to $i(v)$, and an interior point of a face $f\in S\setminus i(G)$ whose boundary contains $i(e)$. All flags contained in the closure of the face $f$ share the same vertex, ``the interior point''. Hence, by the construction described in the introduction, each map $\mathcal{M}$ induces a triangulation of the surface $S$. From a combinatoric point of view, one can identify each flag $\Phi$ of the map $\mathcal{M}$ with an \textbf{ordered incident triplet} conformed by a vertex, an edge and a face mutually incident in the map $\mathcal{M}$, it means $\Phi=(i(v),i(e),f)$, where the vertex $i(v)$ is inciden with the edge $i(e)$, which is belong to the boundary of the face $f$. To each flag $\Phi$ of the map $\mathcal{M}$, there exists a unique adjacent flag $\Phi^0$  of the map $\mathcal{M}$ that differs from $\Phi$ only on the vertex, and in the same manner, there exist unique adjacent flags $\Phi^1$ and $\Phi^2$  that differ from $\Phi$ on the edge and on the face, respectively. The flag $\Phi^j$ will be called the \textbf{$j$-adjacent flag of} $\Phi$, with $j\in\{0,1,2\}$. We shall denote by $\mathcal{F(M)}$ the set conformed by all flags of the map $\mathcal{M}$. In Figure \ref{Im:flags}, we show an example of a map on the torus with some flags marked with an arbitrary base flag $\Phi$ and its three $i$-adjacent flags.

%%%%%%%%%%%%%%%%%%%%%%%%%%%%%%%%%%%%%%%%%%%%%%%%%%%%%%%%%%

\subsection{Automorphism and monodromy groups}

An \textbf{automorphism} $h$ of a map $\mathcal{M}$ is an automorphism of the graph $G$, such that it can be extended to a homeomorphism $\widehat{h}$ of the surface $S$ to itself, this is $i \circ h=\widehat{h}\circ i$. The automorphism set of a map $\mathcal{M}$, which will be denoted by $Aut(\mathcal{M})$ has a group structure with the composition operation. Hence, the automorphism group of the map $\mathcal{M}$ is a subgroup of the group of automorphism of the graph $G$, $Aut(\mathcal{M})\leq Aut(G)$.  The automorphism group $Aut(\mathcal{M})$ acts on the set of flags $\mathcal{F}(\mathcal{M})$. Namely, this  action is \emph{free}; that is, each element of $Aut (\mathcal{M})$ is completely determined by the image of a given flag (see \cite[Lemma 3.1]{Graver}). Hence, $\mathcal{O}_{\Phi}$ will denote the orbit of each flag $\Phi \in \mathcal{F}(\mathcal{M})$ under the action of the automorphism group $Aut(\mathcal{M})$, and we denote by 
\begin{equation}\label{action}
\mathcal{O}rb(\mathcal{M}):=\{\mathcal{O}_{\Phi}\mid\Phi \in \mathcal{F}(\mathcal{M})\}
\end{equation}
the set conformed by the orbits defined by the action of $Aut(\mathcal{M})$ on $\mathcal{F(M)}$.

A map $\mathcal{M}$ is called \textbf{$k$-orbit map} if the action of its automorphism group $Aut(\mathcal{M})$ induces $k$ orbits on the set of flags $\mathcal{F}(\mathcal{M})$, for some $k\in\mathbb{N}$ (see \cite[Section 3]{Pellicer}). In the literature, a map $\mathcal{M}$ is called \textbf{regular}, if the action of $Aut(\mathcal{M})$ on $\mathcal{F(M)}$ induces one orbit on the set of flags. And a map $\mathcal{M}$ is called \textbf{chiral},  if the action of $Aut(\mathcal{M})$ on $\mathcal{F(M)}$ induces two orbits on the set of flags, with the property that all adjacent flags belong to different orbits (see \emph{e.g.},  \cite{McSc}).

We denote as $s_j: \mathcal{F(M)}  \to \mathcal{F(M)}$, for every $j\in\{0,1,2\}$, the permutation on the set of flags $\mathcal{F(M)}$ of the map $\mathcal{M}$, which sends each flag $\Phi$ to its $j$-adjacent flag $\Phi^j$,
\begin{equation}\label{j_adjacent}
\Phi \to \Phi \cdot s_j:=\Phi^j.
\end{equation}

The permutation $s_j$ is an involution, it means 
\[
(\Phi\cdot s_j) \cdot s_j=\Phi^j\cdot s_j=(\Phi^j)^j=\Phi, \text{ for each flag in } \mathcal{M}.
\]
Moreover, $s_j$ is not an automorphism of the map $\mathcal{M}$ because it does not induce a homeomorphism of the surface $S$, it is merely a bijection in the set of flags (see \emph{e.g.} \cite[Section 2]{CPRRV}).

The \textbf{monodromy group}\footnote{There are some other authors that now prefer to refer to this group as the connection group. Stephen E. Wilson was the one introducing the subject like this.} $Mon(\mathcal{M})$ of the map $\mathcal{M}$ is the subgroup of the permutation group of  the set of flags $\mathcal{F(M)}$, which is generated by the elements $s_0$, $s_1$ and $s_2$, \emph{i.e.}, 
\begin{equation}\label{mon}
Mon(\mathcal{M}):=\langle s_0,s_1,s_2\rangle,
\end{equation}

Let $\Phi$ be a flag of the map $\mathcal{M}$ and let $j_0$ and $j_{1}$ be index in the set $\{0,1,2\}$, then by equation (\ref{j_adjacent}) we introduce the following notation
\begin{equation}
(\Phi\cdot s_{j_0})\cdot s_{j_1}=(\Phi^{j_0})\cdot s_{j_1}:=\Phi^{j_0,j_1}.
\end{equation}
Hence, one can naturally define the right action of $Mon(\mathcal{M})$ on $\mathcal{F(M)}$ as follows 
\begin{equation}\label{monodromyaction}
\alpha(w, \Phi):=\Phi\cdot w,
\end{equation}
for each $\Phi\in\mathcal{F(M)}$ and each  $w\in Mon(\mathcal{M})$. Where, for each $w\in Mon(\mathcal{M})$ there are integers $j_{0},j_{1},\ldots,j_{k} \in \{0,1,2\}$, for any $k \in \mathbb{N}$, such that $w=s_{j_{0}}\circ s_{j_{1}}\circ \ldots \circ s_{j_{k}}$, then the equation (\ref{monodromyaction}) can be written as
\[
\Phi\cdot w=\Phi \cdot ( s_{j_{0}}\circ s_{j_{1}}\circ \ldots \circ s_{j_{k}})=\Phi^{j_{0},j_{1},\ldots,j_{k}}.
\]
In fact, this group satisfies the following properties (see \emph{e.g.}, \cite{Maniplexes}, \cite{HuOrIv}).
\begin{enumerate}
           \item Its only defined relations are $s_i^2=Id$ for each $i\in\{0,1,2\}$ and, $(s_i\circ s_j)^2=Id$ whenever $|i-j|\geq 2$ such that $i,j\in\{0,1,2\}$. In other words, The elements $s_0, s_1, s_2$ and $s_0\circ s_2$ are fixed-point free involutions.
           \item The group $Mon(\mathcal{M})$ acts trasitively on $\mathcal{F(M)}$.
	%\item The group $Mon(\mathcal{M})$ is transitive on $\mathcal{F(M)}$.
	%\item The elements $s_0, s_1$ and $s_2$ are fixed-point free involutions.
	%\item $s_0\circ s_2=s_2\circ s_0$, and $s_0\circ s_2$ is fixed-point free involution. 
\end{enumerate}

%%%%%%%%%%%%%%%%%%%%%%%%%%%%%%%%%%%%%%%%%%%%%%%%%%%%%%%%%%%%%%
%%%%%%%%%%%%%%%%%%%%%%%%%%%%%%%%%%%%%%%%%%%%%%%%%%%%%%%%%%%%%%
%%%%%%%%%%%%%%%%%%%%%%%%%%%%%%%%%%%%%%%%%%%%%%%%%%%%%%%%%%%%%%

\section{Pregraph and Symmetry type graph}\label{sec:pregraph}

\begin{figure}[h!]
\centering
\begin{tabular}{ccc}
\includegraphics[scale=0.4]{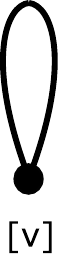}& &\includegraphics[scale=0.4]{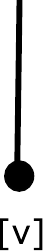}\\
a. Loop of the  && b. Semi-edge of the \\
  pregraph. & & pregraph.\\
  \end{tabular}
 \caption{\emph{An edge of the pregraph $G_{\mathcal{B}}$ with the same vertices.}} 
  \label{semi_edge}
\end{figure}

Given a graph $G$, we consider an edge colouring $\mathcal{C}$ of $\mathrm{G}$ and a partition $\mathcal{B}$ of the vertex set $\mathrm{V}(\mathrm{G})$ of the graph. The coloured quotient with respect to $\mathcal{B}, G_{\mathcal{B}},$ is defined as the pregraph with vertex set $\mathcal{B},$ such that for any two vertices $B, C \in \mathcal{B},$ there is an edge of colour $a$ from $B$ to $C$ if and only if there exist two classes  $[u],[v] \in\mathcal{B}$ define an edge with colour $k$, if and only if, there exist $\hat{u} \in [u]$ and $\hat{v}\in [v]$ such that there is an edge with colour $k$ from $\hat{u}$ to $\hat{v}$. It could be that an edge with colour $k$ of the pregraph $G_{\mathcal{B}}$ has the same vertices \emph{i.e.}, $[u]=[v]$, in this case the edge ``loop'' will be called \textbf{semi-edge} with colour $k$ and it will be thought as is shown in Figure \ref{semi_edge}. For more details, we refer the reader to \cite[Section 3]{Maniplexes}.
We will denote an edge of $G_{\mathcal{B}}$ having colour $k$, and vertices the class $[v]$ and $[u]$ as $\left\{[u],[v]\right\}_k$. Similarly, we will denote an semi-edge of $G_{\mathcal{B}}$ having colour $k$, and vertices the class $[v]$ as $\left\{[v]\right\}_k$.

Given a map $\mathcal{M}$, then the \textbf{flag graph} ${G}_{\mathcal{M}}$ \textbf{corresponding to} $\mathcal{M}$ is the graph whose  set of vertices is conformed by the flags of the map $\mathcal{M}$, and two flags $\Phi,\Psi\in \mathcal{F(M)}$ define an edge if they are adjacent. The flag graph $G_{\mathcal{M}}$ is $3$-regular \emph{i.e.}, each vertex of $G_{\mathcal{M}}$ has degree three, because each flag $\Phi$ of $\mathcal{M}$ is only adjacent to three flags: $\Phi^0$, $\Phi^1$ and $\Phi^2$. Thus, we consider the three different colours $k_1, k_2$ and $k_3$, and define the edge colouring of $G_{\mathcal{M}}$
\begin{equation}\label{edge_colouring}
\mathcal{C}  :  E(G_{\mathcal{M}})  \to \left\{k_0,k_1,k_2\right\},
\end{equation}
which sends the edge with vertices on the flag $\Phi,\Psi$ to the colour $k_0$, $k_1$ or $k_2$, if they differ by a vertex, an edge or a face, respectively.

The function $\mathcal{C}$ sends the edges of $G_{\mathcal{M}}$ with vertices $\Phi$ and $\Phi^j$ to the colour $k_j$, with $j\in\{0,1,2\}$, for each $\Phi\in\mathcal{F(M)}$. Moreover, the colour preserving automorphism group of $G_{\mathcal{M}}$ is isomorphic to the automorphism group of $\mathcal{M}$ (see \cite[Subsection 1.4]{Gunnar}).

Given a map $\mathcal{M}$, let $\mathcal{C}$ be the edge colouring of the flag graph $G_{\mathcal{M}}$ defined in equation (\ref{edge_colouring}), and $\mathcal{O}rb(\mathcal{M})$ the set of orbits defined in equation (\ref{action}). Then \textbf{the symmetry type graph $\mathcal{T(M)}$ associated to $\mathcal{M}$} is the pregraph of $G_{\mathcal{M}}$ with respect to $\mathcal{O}rb(\mathcal{M})$. The graph $\mathcal{T(M)}$ is such that its set of vertices is $\mathcal{O}rb(\mathcal{M})$, and two orbits  $\mathcal{O}_{\Phi},\mathcal{O}_{\Psi}\in \mathcal{O}rb(\mathcal{M})$ define an edge with colour $k_j$, with $j\in\{0,1,2\}$, if and only if, there exist $\hat{\Phi} \in \mathcal{O}_{\Phi}$ and $\hat{\Psi}\in \mathcal{O}_{\Psi}$ such that there is an edge with colour $k_j$ from $\hat{\Phi}$ to $\hat{\Psi}$. We note that if $\mathcal{M}$ is $k$-orbit map, then $\mathcal{T(M)}$ has exactly $k$ vertices.

%For example, the symmetry type graph $\mathcal{T(M)}$ of a regular map $\mathcal{M}$ has a single vertex $\mathcal{O}_{\Phi}$ and three semi-edges $\{\mathcal{O}_{\Phi}\}_{j}$ having colours $k_{j}$ respectively, for each  $j=\{0,1,2\}$, because the automorphism $Aut(\mathcal{M})$ defined only one orbit $\mathcal{O}_{\Phi}$, and for any flag $\Phi$ of the map $\mathcal{M}$ and its respective $j$-adjacent flag $\Phi^j$ belong to the orbit $\mathcal{O}_{\Phi}$ (see Figure \ref{pregraph}-a). On the other hand, the symmetry type graph $\mathcal{T(M)}$ of a chiral map $\mathcal{M}$ has exactly two vertices $\mathcal{O}_{\Phi}$ and $\mathcal{O}_{\Psi}$ joined by three edges $\{\mathcal{O}_{\Phi},\mathcal{O}_{\Psi}\}_{j}$ having colours $k_{j}$ respectively, for each $j\in\{0,1,2\}$, because the automorphism $Aut(\mathcal{M})$ defined two orbits  $\mathcal{O}_{\Phi}$ and $\mathcal{O}_{\Psi}$ such that two adjacent flags belong to different orbits (see Figure \ref{pregraph}-b).
%
%\begin{figure}[H]
%\centering
%\begin{tabular}{ccc}
%\includegraphics[scale=0.35]{pregraph_regular}& &\includegraphics[scale=0.35]{pregraph_quiral}\\
%a. Symmetry type graph of a && b. Symmetry type graph of a\\
%  regular map. & &  chiral map.\\
%  \end{tabular}
% \caption{\emph{Classic examples of symmetry type graphs.}} 
%  \label{pregraph}
%\end{figure}

\begin{figure}[h!]
  \centering
  \includegraphics[scale=0.45]{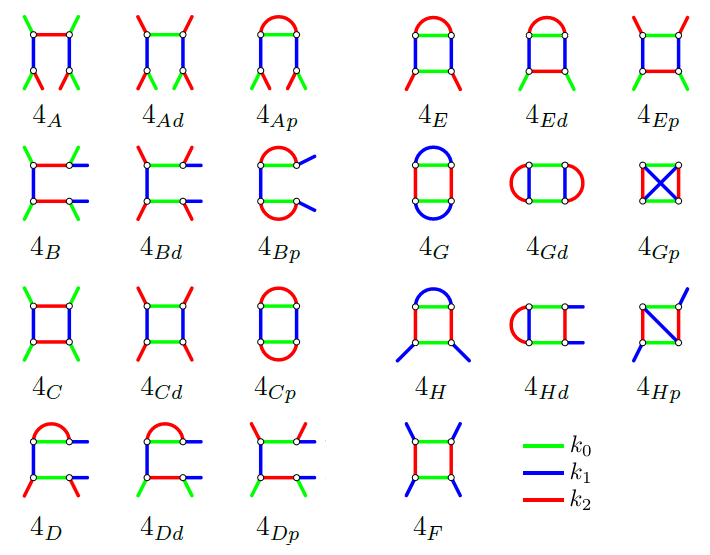}
  \caption{\emph{Symmetry type graphs associated to the $4$-orbit maps, with edges and semi-edges with colours $k_{j}$ for $j\in\{0,1,2\}$.}}
  \label{classes}
\end{figure}

There are twenty two symmetry type graphs associated to $4$-orbit maps (see \cite{Pellicer}). It means, the symmetry type graph of any $4$-orbit map is isomorphic to one of those pregraph shown in Figure \ref{classes}.

If $\mathcal{M}$ is a map and $\mathcal{C}$ is the edge colouring of the flag graph $G_{\mathcal{M}}$, then the set of edges of $G_{\mathcal{M}}$ with colour $j$ forms a perfect matching, for $j\in\{0,1,2\}$ (see \cite[Section 2]{Hubard4}). Hence, the graph $G_{\mathcal{M}}^{j,i}$ conformed by the edges of $G_{\mathcal{M}}$ with colours $j$ and $i$, such that $j\neq i\in\{0,1,2\}$, is a subgraph of $G_{\mathcal{M}}$  whose connected components are even cycles. The graph $G_{\mathcal{M}}^{j,i}$ is called a \textbf{$2$-factor of $G_{\mathcal{M}}$}. 

\begin{figure}[h]
	\centering
	\includegraphics[scale=0.25]{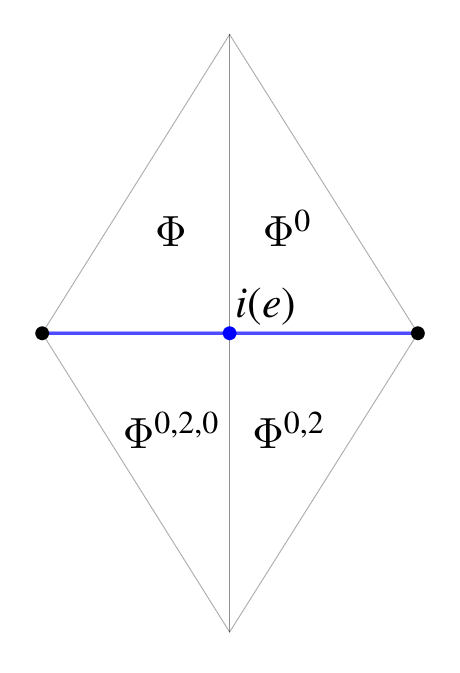}
	\caption{\emph{Orbit $\mathcal{O}_{\Phi}^{\langle s_{0},s_{2} \rangle}$ around the edge $i(e)$.}}	
	\label{sequenceE1}
\end{figure}

Given that the permutation $s_{0}\circ s_{2}$ of the monodromy group $Mon(\mathcal{M})$ is fixed-point free involution, then the cycles of the subgraph $G_{\mathcal{M}}^{0,2}\subset G_{\mathcal{M}}$ are colourable alternantely with colours $k_0$ and $k_2$, and all them have length four. Hence, if $\Phi$ is a flag of $\mathcal{M}$, then the elements in the sequence of flags $\Phi, \Phi^{0},\Phi^{0,2},\Phi^{0,2,0}$ are the vertices of a cycle of $G_{\mathcal{M}}^{0,2}$, where the second coordinate of the flags $\Phi, \Phi^{0},\Phi^{0,2},\Phi^{0,2,0}$ are the same. Therefore, there is a one-to-one correspondence between  the set of edges of $\mathcal{M}$ and the cycles of $G_{\mathcal{M}}^{0,2}$. This correspondence is given by the orbits of $\mathcal{F(M)}$ under the action of the subgroup of $Mon(\mathcal{M})$ generated by $s_0$ and $s_1$, \emph{i.e.},
\begin{equation}\label{eq_orbit_edge}
i(e) \to \{\Phi, \Phi^{0},\Phi^{0,2},\Phi^{0,2,0}\}:=\mathcal{O}_{\Phi}^{\langle s_{0},s_{2} \rangle}=\{\Phi\cdot w : w \in \langle s_{0},s_{2} \rangle\}, 
\end{equation}
being $\Phi$ a flag of $\mathcal{M}$, such that its second coordinate is $i(e)$. We will say that \textbf{the orbit $\mathcal{O}_{\Phi}^{\langle s_{0},s_{2} \rangle}$ is around the edge $i(e)$} (see Figure \ref{sequenceE1}). Therefore, the cycle of $G_{\mathcal{M}}^{0,2}$ such that its vertices are the flags belong to $\mathcal{O}_{\Phi}^{\langle s_{0},s_{2} \rangle}$ will be denoted as $C_{i(e)}$.

\begin{figure}[H]
	\centering
	\includegraphics[scale=0.25]{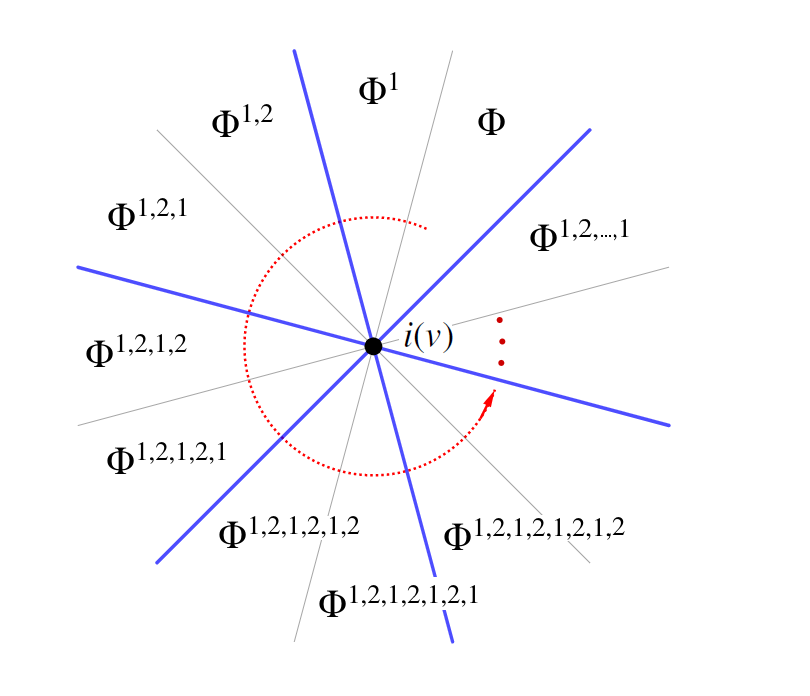}
	\caption{\emph{Orbit $\mathcal{O}_{\Phi}^{\langle s_{1},s_{2} \rangle}$ around the vertex $i(v)$.}}
	\label{sequenceV1}	
\end{figure}

Analogously, the permutation $s_{1}\circ s_{2}$ of the monodromy group $Mon(\mathcal{M})$ is fixed-point free, and it has finite order, then the cycles of the subgraph $G_{\mathcal{M}}^{1,2}\subset G_{\mathcal{M}}$ are colourable alternantely with colours $k_1$ and $k_2$, and all them have even length.  Hence, if $\Phi$ is a flag of $\mathcal{M}$, then the elements in the finite sequence of flags $\Phi, \Phi^{1},\Phi^{1,2},\ldots, \Phi^{1,2,\ldots,1}$ are the vertices of a cycle of $G_{\mathcal{M}}^{1,2}$, where the first coordinate of the flags $\Phi, \Phi^{1},\Phi^{1,2},\ldots, \Phi^{1,2,\ldots,1}$ are the same. Therefore, there is a biunique correspondence between  the set of vertices of $\mathcal{M}$ and the cycles of $G_{\mathcal{M}}^{1,2}$. This correspondence  is given by the orbits of $\mathcal{F(M)}$ under the action of the subgroup of $Mon(\mathcal{M})$ generated by $s_1$ and $s_2$, \emph{i.e.},
\begin{equation}\label{eq_orbit_vertex}
i(v) \to \{\Phi, \Phi^{1},\Phi^{1,2},\ldots, \Phi^{1,2,\ldots,1}\}:=\mathcal{O}_{\Phi}^{\langle s_{1},s_{2} \rangle}=\{\Phi\cdot w : w \in \langle s_{1},s_{2} \rangle\}, 
\end{equation}
being $\Phi$ a flag of $\mathcal{M}$, such that its first coordinate is $i(v)$. We will say that \textbf{the orbit $\mathcal{O}_{\Phi}^{\langle s_{1},s_{2} \rangle}$ is around the vertex $i(v)$} (see Figure \ref{sequenceV1}). Therefore, the cycle of $G_{\mathcal{M}}^{1,2}$ such that its vertices are the flags belong to $\mathcal{O}_{\Phi}^{\langle s_{1},s_{2} \rangle}$ will be denoted as $C_{i(v)}$.

\begin{figure}[H]
	\centering
	\includegraphics[scale=0.25]{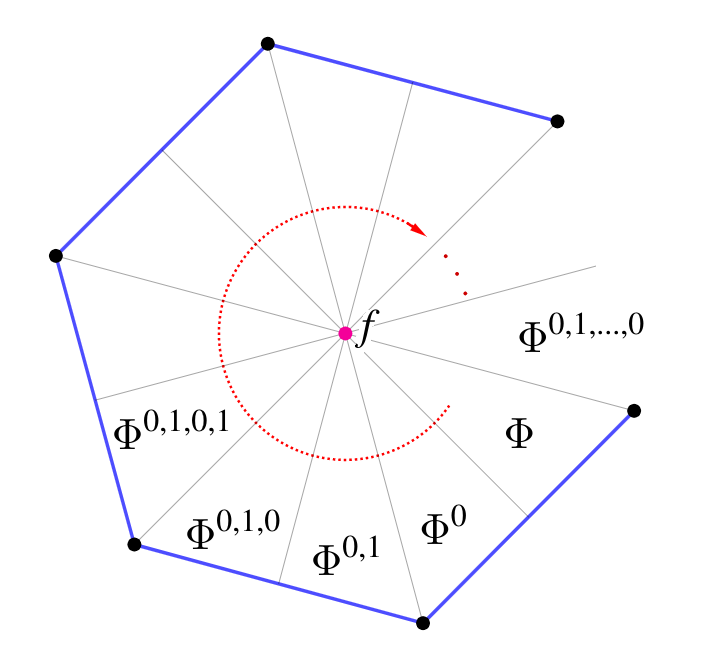}
	\caption{\emph{Orbit $\mathcal{O}_{\Phi}^{\langle s_{0},s_{1} \rangle}$ around the face $f$.}}	
	\label{sequenceF1}
\end{figure}

Likewise, the permutation $s_{0}\circ s_{1}$ of the monodromy group $Mon(\mathcal{M})$ is fixed-point free, and it has finite order, then the cycles of the subgraph $G_{\mathcal{M}}^{0,1}\subset G_{\mathcal{M}}$ are colourable alternantely with colours $k_0$ and $k_1$, and all them have even length.  Hence, if $\Phi$ is a flag of $\mathcal{M}$, then the elements in the finite sequence of flags $\Phi, \Phi^{0},\Phi^{0,1},\ldots, \Phi^{0,1,\ldots,0}$ are the vertices of a cycle of $G_{\mathcal{M}}^{0,1}$, where the third coordinate of the flags $\Phi, \Phi^{0},\Phi^{0,1},\ldots, \Phi^{0,1,\ldots,0}$ are the same. Therefore, there is a biunique correspondence between  the set of faces of $\mathcal{M}$ and the cycles of $G_{\mathcal{M}}^{0,1}$. This correspondence is given by the orbits of $\mathcal{F(M)}$ under the action of the subgroup of $Mon(\mathcal{M})$ generated by $s_0$ and $s_1$, \emph{i.e.},
\begin{equation}\label{eq_orbit_face}
f \to \{\Phi, \Phi^{0},\Phi^{0,1},\ldots, \Phi^{0,1,\ldots,0}\}:=\mathcal{O}_{\Phi}^{\langle s_{0},s_{1} \rangle}=\{\Phi\cdot w : w \in \langle s_{0},s_{1} \rangle\}, 
\end{equation}
being $\Phi$ a flag of $\mathcal{M}$, such that its third coordinate is $f$. We will say that \textbf{the orbit $\mathcal{O}_{\Phi}^{\langle s_{0},s_{1} \rangle}$ is around the face $f$}, (see Figure \ref{sequenceF1}). Therefore, the cycle of $G_{\mathcal{M}}^{0,1}$ such that its vertices are the flags belong to $\mathcal{O}_{\Phi}^{\langle s_{0},s_{1} \rangle}$ will be denoted as $C_{f}$.

\subsection{Dual and petrie-dual maps}

Let $\mathcal{M}$ and $\mathcal{N}$ be two maps and let $\mathcal{F(M)}$ and $\mathcal{F(N)}$ be their respectively set of flags, 
a \textbf{duality} $\delta$ from $\mathcal{M}$ to $\mathcal{N}$ is a bijection funtion $\delta: \mathcal{F(M)}\to\mathcal{F(N)}$,  satisfying $\Phi^{i}\delta= (\Phi \delta)^{i}$ for each flag $\Phi$ of  $\mathcal{F(M)}$ and each $i\in\{0,1,2\}$. The map $\mathcal{N}$ is called the \textbf{dual map} of the map $\mathcal{M}$, if there is a duality from $\mathcal{M}$ to $\mathcal{N}$, and we shall denote it as $\mathcal{M}^{*}$. If there exists a duality from the map $\mathcal{M}$ to itself, then the map $\mathcal{M}$ will called \textbf{self-dual}. Note that the duality $\delta$ defines a bijection from the vertices of the symmetry type graph $\mathcal{T(M)}$ to the vertices of the symmetry type graph $\mathcal{T(M^{*})}$, which sends the edge of color $i$ of $\mathcal{T(M)}$ onto the edge of color $2-i$ of $\mathcal{T(M^{*})}$, for each $i\in\{0,1,2\}$.

A \textbf{Petrie polygon} in a map $\mathcal{M}$ is defined as a zig-zag path in the map. More precisely, we
start at a vertex, then go along an edge to an adjacent vertex, the turn left and go to the next
vertex and then turn right, and so on, (or interchange left and right.) We have a path in which
two consecutive edges belong to the same face but no three consecutive edges belong to the
same face \cite{CoMo}. Note that each edge of a Petrie polygon appears either just once in exactly two different Petrie
polygons of the map  $\mathcal{M}$, or twice in the same Petrie polygon of the map $\mathcal{M}$. Hence we can define  a map with the same set of vertices and edges of $\mathcal{M}$, but with the Petrie polygons as faces. This
map is known as the \textbf{Petrie-dual} or \textbf{Petrial} map of $\mathcal{M}$, which will be denoted by $\mathcal{M^P}$. If the map
$\mathcal{M}$ is isomorphic to its respective Petrie-dual map  $\mathcal{M^P}$, then $\mathcal{M}$ is said to be \textbf{self-Petrie}.

We have the following result for the dual and petrie dual maps.

\begin{prop}[\cite{Hubard4}]\label{prop:Hubard13}
 If a map $\mathcal{M}$ has symmetry type graph $\mathcal{T(M)}$, then 

\begin{enumerate}
\item Its dual map $\mathcal{M^{*}}$ has the
dual of $ \mathcal{T(M)}$ as symmetry type graph.

\item Its Petrie-dual $\mathcal{M^P}$ has the petrie-dual of $\mathcal{T(M)}$ as symmetry type graph.
\end{enumerate}
\end{prop}

In the Table  \ref{Table:dual-petrial} is shown the dual and the petrial of each symmetry type graph associated to the $4$-orbit maps (see \cite{Pellicer}).

\begin{table}[h!]
\begin{center}
\begin{tabular}{|c|c|c||| c|c|c| }
\hline 
 Symmetry & Dual &  Petrial &  Symmetry & Dual &  Petrial\\
type graph &        &             &    type graph         &       &  \\
\hline 
$4_A$       & $4_{Ad}$   &  $4_{Ap}$  &  $4_{Ad}$ & $4_A$     & $4_{Ad}$\\
\hline
$4_{Ap}$  & $4_{Ap}$   &  $4_{A} $  &  $4_{B}$   & $4_{Bd}$ & $4_{Bp}$  \\
\hline
 $4_{Bd}$ & $4_{B}$     & $4_{Bd}$  &   $4_{Bp}$  & $4_{Bp}$  & $4_{B}$ \\
\hline
 $4_C$     &  $4_{Cd}$   & $4_{Cp}$ &  $4_{Cd}$   & $4_{C}$    & $4_{Cd}$  \\
\hline
 $4_{Cp}$&   $4_{Cp}$  & $4_{C}$   &  $4_{D}$    & $4_{Dd}$   & $4_{Dp}$  \\
\hline
$4_{Dd}$ &   $4_{D}$    & $4_{Dd}$ &  $4_{Dp}$  & $4_{Dp}$   &  $4_{D}$      \\
\hline
$4_{E}$ &   $4_{Ed}$     & $4_{Ep}$ &  $4_{Ed}$   & $4_{E}$    & $4_{Ed}$   \\
\hline 
$4_{Ep}$&  $4_{Ep}$     & $4_{E}$   &  $4_{F}$     & $4_{F}$    & $4_{F}$    \\
\hline
 $4_{G}$  & $4_{Gd}$    & $4_{Gp}$&  $4_{Gd}$    & $4_{G}$   & $4_{Gd}$   \\
\hline
 $4_{Gp}$ & $4_{Gp}$   & $4_{G}$ &   $4_{H}$      & $4_{Hd}$ & $4_{Hp}$  \\
\hline
$4_{Hd}$  & $4_H$         & $4_{Hd}$ & $4_{Hp}$    &  $4_{Hp}$  & $4_{H}$ \\
\hline
\end{tabular}
\caption{\emph{The dual and the petrial of each symmetry type graph associated to the $4$-orbit maps.}} 
\label{Table:dual-petrial}
\end{center}
\end{table}

%%%%%%%%%%%%%%%%%%%%%%%%%%%%%%%%%%%%%%%%%%%%%%%%%%%%%%%%%%
%%%%%%%%%%%%%%%%%%%%%%%%%%%%%%%%%%%%%%%%%%%%%%%%%%%%%%%%%%
%%%%%%%%%%%%%%%%%%%%%%%%%%%%%%%%%%%%%%%%%%%%%%%%%%%%%%%%%%

\section{Characteristic system of vertices and faces}\label{sec:charsys}

In this section we discuss about the local combinatorial nature of $4$-orbit maps, from the point of view of their symmetry type graph, characteristic system of a vertex and characteristic system of a face. 

%%%%%%%%%%%%%%%%%%%%%%%%%%%%%%%%%%%%%%%%%%%%%%%%%%%%%%%%%%

\subsection{Characteristic system of a vertex}

Consider  a $4$-orbit map $\mathcal{M}$, let $G_{\mathcal{M}}$ be the flag graph associated to $\mathcal{M}$ and $\mathcal{C}$ the edge colouring defined in equation (\ref{edge_colouring}). If $i(v)$ is a vertex of the map $\mathcal{M}$, then there is a cycle $C_{i(v)}$ of $G_{\mathcal{M}}$ around $i(v)$ (see equation (\ref{eq_orbit_vertex})), having even length and being two colourable alternating the colours $k_1$ and $k_2$. From this properties is motivated the following definition.

\begin{definition}\label{m_v}
Let $\mathcal{M}$ be a $4$-orbit map, and $\mathcal{T(M)}$  its symmetry type graph. If  $i(v)$ is a vertex of $\mathcal{M}$,  then the ordered triplet 
\[
(2m_v,k_1,k_2)
\]
associated to $i(v)$ is called the \textbf{characteristic system of the vertex $i(v)$}, where $2m_v$ is the length of the   two colourable alternately cycle $C_{i(v)}$, with colours $k_1$ and $k_2$, for some $m_v\in\mathbb{N}$. The positive integer $m_v$ corresponds to the degree of the vertex $i(v)$.
\end{definition}

\begin{figure}[h!]
  \centering
  \includegraphics[scale=0.45]{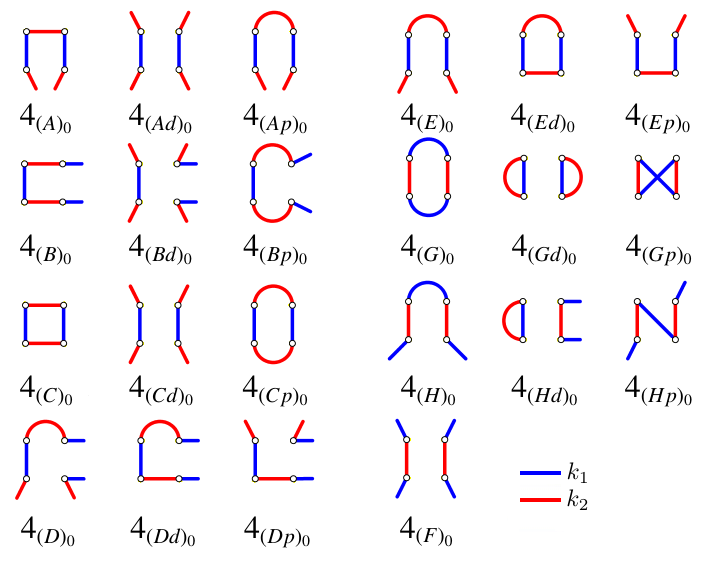}
  \caption{\emph{Pregraphs $\mathcal{T}_0\mathcal{(M)}$ associated to the $4$-orbit maps without the edges and semi-edges with colour $k_0$.}}
  \label{Im:VertexPregraph}
\end{figure}

\begin{figure}[H]
  \centering
  \includegraphics[scale=0.3]{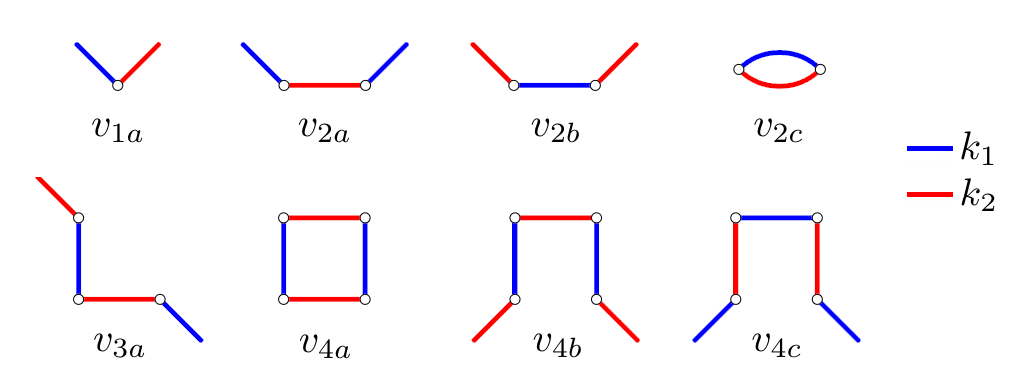}
  \caption{\emph{Pregraphs $v_x$.}}
  \label{Im:VerticesPregraph}
\end{figure}

Given that the characteristic system of the vertex $i(v)$ is determined by three parameters,  if we consider a symmetry type graph $\mathcal{T(M)}$ described in Figure  \ref{classes}, and we remove the edges and semi-edges having colour $k_0$, then we hold a new pregraph $\mathcal{T}_0\mathcal{(M)}$ isomorphic to one of those twenty two pregraphs shown in Figure \ref{Im:VertexPregraph}.
\begin{remark}\label{Cor:No.Vertex}
The pregraph $\mathcal{T}_0\mathcal{(M)}$ is conformed by at most three connected components. Moreover, each connected component of $\mathcal{T}_0\mathcal{(M)}$ is isomorphic to one of those eight pregraphs shown in Figure \ref{Im:VerticesPregraph}, which we denote as  $v_x$, for some $x$ in the set of index $\{1a,2a,2b,2c,3a,4a,4b,4c\}$.
\end{remark}

If we fix a vertex $i(v)$ of the map $\mathcal{M}$, we hold the cycle $C_{i(v)}$ of $G_{\mathcal{M}}$ around $i(v)$, and remember that $C_{i(v)}$ is a two colourable alternately cycle, then we can introduce the set of orbits $\mathcal{O}rb(C_{i(v)})$ defined by the action of $Aut(\mathcal{M})$ on $\mathcal{F(M)}$ restricted to the flags that conformed the vertices of $C_{i(v)}$. Using the definition of a pregraph from Section  \ref{sec:pregraph} we hold that the pregraph $\overline{C}_{i(v)}$ of the cycle $C_{i(v)}$ with respect to $\mathcal{O}rb(C_{i(v)})$ induces the following definition.
 
 \begin{definition}
Consider the pregraph $\overline{C}_{i(v)}$ which is contained into a connected component of $\mathcal{T}_0\mathcal{(M)}$, by construction, $\overline{C}_{i(v)}$ is a pregraph of type $v_x$, for some $x$ in the set of index $\{1a,2a,2b,2c,3a,4a,4b,4c\}$. This connected component is called  \textbf{the vertex type graph $\mathcal{T}(i(v))$ of the vertex $i(v)$}.
\end{definition}

\begin{thm}\label{thm:vertex_cycle}
Let us considerer a $4$-orbit map and let  $i(v)$ be a vertex of the map. If $v_{2a}$ is the vertex type graph of $i(v)$, then the degree of the vertex is even. Moreover, the characteristic system of the vertex is $(4n,k_1,k_2)$ for some $n \geq 2$.  
\end{thm}

\begin{proof}
Let us consider an edge $i(e)$ incident to $i(v)$ and $f$ a face of the map $\mathcal{M}$ such that $i(e)$ belongs to its boundary. We denote as $\Phi$ the flag of the map $\mathcal{M}$ conformed by the triplet 
\[
\Phi:=(i(v),i(e),f).
\]
Suppose that there is a flag $\Psi$ on the map, such that the classes $\mathcal{O}_\Phi, \mathcal{O}_\Psi\in \mathcal{O}rb(\mathcal{M})$ are vertices of the pregraph $v_{2a}$ (see Figure \ref{Vertex_v2a}-a). We will count the number of elements in the set $\mathcal{O}_{\Phi}^{\langle s_{1},s_{2}\rangle}$ (see equation \ref{eq_orbit_vertex}) using the vertex type graph  $v_{2a}$.

 \noindent
Considering the action of $\langle s_{1},s_{2} \rangle$ on the set of flags,  by equation (\ref{eq_orbit_vertex}), it holds that the class 
\[
\mathcal{O}_{\Phi}^{\langle s_{1},s_{2} \rangle}=\{\Phi\cdot w : w\in \langle s_{1},s_{2} \rangle\}
\]
contains all the flags around the vertex $i(v)$, it means that
\[
\mathcal{O}_{\Phi}^{\langle s_{1},s_{2} \rangle}=\{\Phi,\Phi^1, \Phi^{1,2}, \Phi^{1, 2, 1}, \Phi^{1, 2, 1, 2},... \}.
\]

Given that $\{\mathcal{O}_{\Phi}\}_1$ is an edge of the pregraph $v_{2a}$, then the flag $\Phi^1$ belongs to the orbit $\mathcal{O}_{\Phi}$. Without lost of generality we can assume that $\Phi=\Phi_1$ y $\Phi^1=\Phi_2$. Analogously, as the sets $\{\mathcal{O}_{\Phi},\mathcal{O}_{\Psi}\}_2$, $\{\mathcal{O}_{\Psi}\}_1$ and $\{\mathcal{O}_{\Psi},\mathcal{O}_{\Phi}\}_2$ are edges of the pregraph $v_{2a}$, then the flags $\Phi^{1,2}$, $\Phi^{1,2,1}$ and $\Phi^{1,2,1,2}$  belong to the orbit $\mathcal{O}_{\Psi}$, $\mathcal{O}_{\Psi}$ and $\mathcal{O}_{\Phi}$. Hence, we can rewrite  $\Phi^{1,2}=\Psi_1$, $\Phi^{1,2,1}=\Psi_2$ and $\Phi^{1,2,1,2}=\Phi_3$. Following with this construction we obtain the finite sequence $\Phi_1,\Phi_2,\Psi_1,\Psi_2,\Phi_3,...,\Phi_{l-1},\Phi_l,\Psi_{l-1},\Psi_l$ (see Figure \ref{Vertex_v2a}-b), where $\Phi_l \in \mathcal{O}_\Phi$ and $\Psi_l \in \mathcal{O}_\Psi$ for $l \in \{1,2,3,...,n\}$. From this it holds that
\[
\begin{array}{ccl}
\mathcal{O}_{\Phi}^{\langle s_{1},s_{2} \rangle} & = &\{\Phi,\Phi^1,\Phi^{1,2}, \Phi^{1,2,1}, \Phi^{1,2,1,2},... \}\\
&= &\{\Phi_1, \Phi_2, \Psi_1, \Psi_2, \Phi_3,...,\Phi_{l-1},\Phi_l,\Psi_{l-1},\Psi_l\}. \\
\end{array}
\]
\begin{figure}[h]
\centering
\begin{tabular}{ccc}
\includegraphics[scale=0.3]{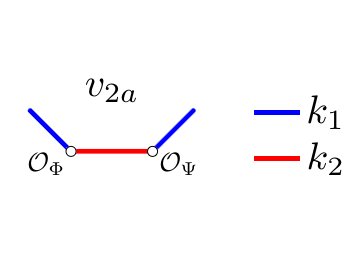}& &\includegraphics[scale=0.3]{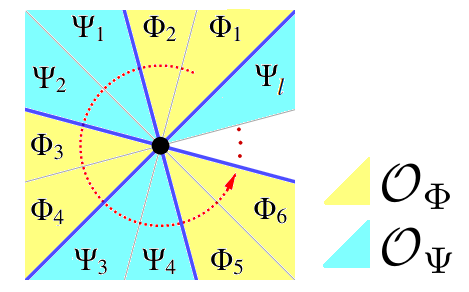}\\
a. Pregraph $v_{2a}$. && b. Orbit $\mathcal{O}_{\Phi}^{\langle s_{1},s_{2} \rangle}$ around\\
      & &  the vertex $i(v)$.\\
  \end{tabular}
  \caption{\emph{Sequence follow by the flags around a vertex $i(v)$ from its vertex type graph of $v_{2a}$.}} 
  \label{Vertex_v2a} 
\end{figure}

 If $\Phi \in \mathcal{O}_{\Phi}^{\langle s_{1},s_{2} \rangle}$, as the pregraph $v_{2a}$ has a semi-edge in $\mathcal{O}_{\Phi}$, then the $1$-adjacent flag $\Phi^{1} \in \mathcal{O}_{\Phi}^{\langle s_{1},s_{2} \rangle}$. Now, we define the following equivalence relation $\sim$ in  $\mathcal{O}_{\Phi}^{\langle s_{1},s_{2} \rangle}$, the flags $\Gamma$ and $\Delta$ of  $\mathcal{O}_{\Phi}^{\langle s_{1},s_{2} \rangle}$ are equivalent if they are $1$-adjacent flags. We remark that each equivalent class $[\Gamma]$ of the quotient set $\mathcal{O}_{\Phi}^{\langle s_{1},s_{2} \rangle}/\sim$ is conformed by exactly two flags of $\mathcal{O}_{\Phi}^{\langle s_{1},s_{2} \rangle}$. This fact implies that the numbers of flags in $\mathcal{O}_{\Phi}^{\langle s_{1},s_{2} \rangle}$ is twice the number of equivalent classes in the quotient set  $\mathcal{O}_{\Phi}^{\langle s_{1},s_{2} \rangle}/\sim$, it means  ${\rm card}\left(\mathcal{O}_{\Phi}^{\langle s_{1},s_{2} \rangle}\right)=2{\rm card}\left(\mathcal{O}_{\Phi}^{\langle s_{1},s_{2} \rangle}/\sim\right)$. Given that  there are at least four  flags $\Phi_1,\Phi_2,\Psi_1,\Psi_2$ in $\mathcal{O}_{\Phi}^{\langle s_{1},s_{2} \rangle}$, then by definition \ref{m_v} it follows that $m_v=2n$  for any positive integer $n\geq 2$. Thus, we conclude that  the characteristic system of the vertex $i(v)$ is $(4n, k_1, k_2)$.

%So we can say that in $\mathcal{O}_{\Phi}^{\langle s_{1},s_{2} \rangle}$ two elements $\phi, \beta$ are equivalent if and only if {\color{red}$\Phi \cap \beta \Leftrightarrow \beta = \Phi^{1}$, that is, $\beta$ is the 1-adjacent flag of $\Phi$, andunder this relation the class $[\Phi]$ has only two elements in $\mathcal{O}_{\Phi}^{\langle s_{1},s_{2} \rangle}$. 

%Hence the cardinalyty of the orbit and the generator of the class are related by
%
%$$ |\mathcal{O}_{\Phi}^{\langle s_{1},s_{2} \rangle}|=2 |\mathcal{O}_{\Phi}^{\langle s_{1},s_{2} \rangle}/\sim|.$$}

%Furthemore, the smallest sequence we can get is $\Phi_1,\Phi_2,\Psi_1,\Psi_2$, then
%
 %$$ |\mathcal{O}_{\Phi}^{\langle s_{1},s_{2} \rangle}|=2 |\mathcal{O}_{\Phi}^{\langle s_{1},s_{2} \rangle}/\sim|>4,$$ 
 %
% it means $m_v=2n$ for some natural number $n\geq 2$. From this, it holds that the characteristic system of the vertex is $(4n, k_1, k_2)$.
\end{proof}

\begin{figure}[h!]
  \centering
  \includegraphics[scale=0.28]{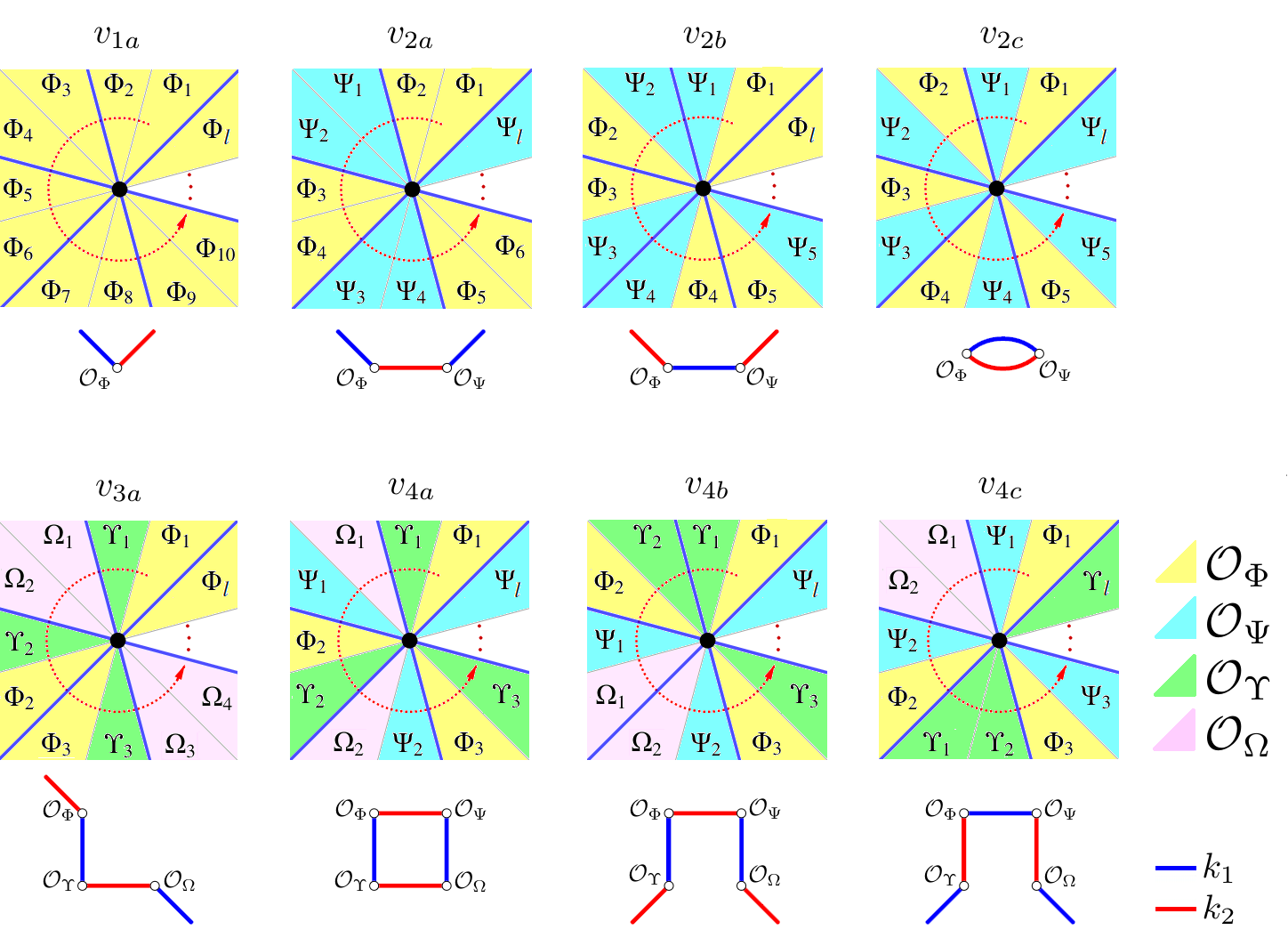}
  \caption{\emph{Sequence follow by the flags around a vertex $i(v)$ from its vertex type graph.}}
  \label{Im:Vertices}
\end{figure}

If we consider any vertex $i(v)$ of the $4$-orbit map $\mathcal{M}$ and we suppose that its vertex type graph $\mathcal{T}(i(v))$ associated is $v_x$, for any $x$ in the set of index $\{1a, 2a, 2b, 2c, 3a, 4a, 4b, 4c\}$, then following the same ideas in the proof of the Theorem \ref{thm:vertex_cycle} it is easy to find the degree of $i(v)$ and the characteristic system of $i(v)$. These results are collected in  Table \ref{Table:VertexType}. In Figure \ref{Im:Vertices} are represented the eight different vertex type graphs. 

\begin{table}[h]
\begin{center}
\begin{tabular}{|c|c|c|c|}
\hline 
 Vertex type graph &  Degree & & Characteristic System \\
of the vertex $i(v)$&   of the vertex $i(v)$ & &of the vertex $i(v)$\\ 
\hline 
$v_{1a}$ & $n$ &$n \geq 3$& $(2n,k_1,k_2)$  \\
\hline 
$v_{2a}$ & $2n$ &$n \geq 2$ & $(4n,k_1,k_2)$  \\
\hline 
$v_{2b}$ & $2n$ &$n \geq 2$& $(4n,k_1,k_2)$  \\
\hline 
$v_{2c}$ & $n$ &$n\geq 3$& $(2n,k_1,k_2)$  \\
\hline
$v_{3a}$ & $3n$ &$n \geq 1$& $(6n,k_1,k_2)$  \\
\hline 
$v_{4a}$ & $2n$ &$n \geq 2$& $(4n,k_1,k_2)$  \\
\hline 
$v_{4b}$ & $4n$ &$n \geq 1$& $(8n,k_1,k_2)$  \\
\hline 
$v_{4c}$ & $4n$ &$n \geq 1$ & $(8n,k_1,k_2)$  \\
\hline 
\end{tabular}
\caption{\emph{Degree and characteristic system of a vertex  $i(v)$ from its vertex type graph.}} 
\label{Table:VertexType}
\end{center}
\end{table} 

From the Table \ref{Table:VertexType} we hold the following corollary.

\begin{corolario}\label{cor:tipo_impar}
	If a $4$-orbit map $\mathcal{M}$ has a vertex of odd degree, then its symmetry type graph $\mathcal{T(M)}$ is either: $4_{Bd}$, $4_D$, $4_{Dp}$, $4_{Gd}$, or $4_{Hd}$.
\end{corolario}

From the Proposition \ref{prop:Hubard13} and Table \ref{Table:dual-petrial} follows that

\begin{corolario}
If $\mathcal{M}$ is one $4$-orbit map with symmetry type graph $\mathcal{T(M)}$ as given in corolally \ref{cor:tipo_impar}, then 

\begin{enumerate}
\item Its dual map $\mathcal{M^{*}}$ has symmetry type graph either: $4_B$, $4_{Dd}$, $4_{Dp}$, $4_{G}$, or $4_{H}$, respectively.

\item Its Petrie-dual $\mathcal{M^P}$ has symmetry type graph either: $4_{Bd}$, $4_{Dp}$, $4_{D}$, $4_{Gd}$, or $4_{Hd}$, respectively.
\end{enumerate}

\end{corolario}

\begin{remark}\label{Re:VertexType}
Let $i(v_1), i(v_2)$ be vertices of a $4$-orbit map $\mathcal{M}$, and let $C_{i(v_1)},C_{i(v_2)}$ be the cycles of the flag graph $G_{\mathcal{M}}$ associated to the vertices $i(v_1)$ and $i(v_2)$, respectively. Suppose that the pregraphs $\overline{C}_{i(v_1)}$ and $\overline{C}_{i(v_2)}$, are contained into some connected components of $\mathcal{T}_{0}(\mathcal{M})$. Given that the automorphism group $Aut(\mathcal{M})$ acts freely on the set of flags  $\mathcal{F(M)}$, then the vertices $i(v_1)$ and $i(v_2)$ have the same degree, and their characteristic systems are the same
\[
(2m_{v_1},k_1,k_2)=(2m_{v_2},k_1,k_2).
\]
However, if the pregraphs $\overline{C}_{i(v_1)}$ and $\overline{C}_{i(v_2)}$  belong to different connected component of $\mathcal{T}_0(\mathcal{M})$ but they are isomorphic, then the degree of the vertices $i(v_1)$ and $i(v_2)$ are multiples of $s$, for any $s \in\mathbb{N}$.
If $l$ is the number of connected components of $\mathcal{T}_0(\mathcal{M})$, for any $l\in\{1,2,3\}$, then there are $l$ values to the degree of the vertices of $\mathcal{M}$. Hence, there are $n_i$ positive integers,  with $i\in\{1,\ldots,l\}$ such that the characteristic system of any vertex of the map is either $(2n_1,k_1,k_2),\ldots,(2n_l,k_1,k_2)$.
\end{remark}

%%%%%%%%%%%%%%%%%%%%%%%%%%%%%%%%%%%%%%%%%%%%%%%%%%%%%%%%%%%%%%%%%%%
%%%%%%%%%%%%%%%%%%%%%%%%%%%%%%%%%%%%%%%%%%%%%%%%%%%%%%%%%%%%%%%%%%%%
%%%%%%%%%%%%%%%%%%%%%%%%%%%%%%%%%%%%%%%%%%%%%%%%%%%%%%%%%%%%%%%%%%%%

\subsection{Characteristic system of a face}

Let $f$ be a face of the $4$-orbit map, let $G_{\mathcal{M}}$ be the flag graph associated to $\mathcal{M}$ and let $\mathcal{C}$ be the edge colouring defined in equation (\ref{eq_orbit_face}). If $f$ is a face of the map $\mathcal{M}$, then there is a cycle $C_{f}$ of $G_{\mathcal{M}}$ around $f$ (see equation  (\ref{eq_orbit_face})), having length even and being two  colourable alternately with colours $k_0$ and $k_1$. From this properties is motivated the following definition.

\begin{definition}
Let $\mathcal{M}$ be a $4$-orbit map, and $\mathcal{T(M)}$ its symmetry type graph. If $f$ is a face of $\mathcal{M}$, then the ordered triplet
\[
(2m_f,k_0,k_1)
\] 
associated to $f$ is called the \textbf{characteristic system of the face $f$}, where $2m_f$ is the length of the two colourable alternately cycle $C_f$, with colours $k_0$ and $k_1$, for some $m_f \in \mathbb{N}$.
\end{definition}

The positive integer $m_f$ corresponds to the number of edges of the map $\mathcal{M}$ in the boundary of the face $f$. This number will be called \textbf{the size of the face $f$}.

\begin{figure}[h!]
  \centering
  \includegraphics[scale=0.45]{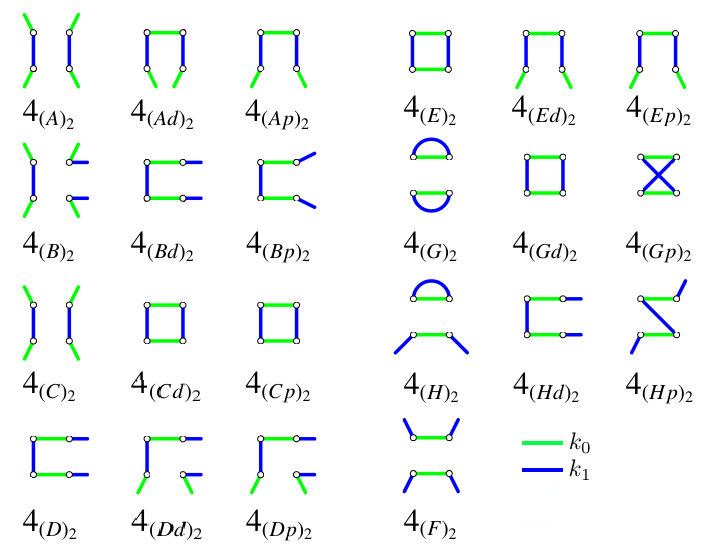}
  \caption{\emph{Pregraphs $\mathcal{T}_2\mathcal{(M)}$ associated to the $4$-orbit maps without the edges and semi-edges with colour $k_2$.}}
  \label{Im:FacePregraph}
\end{figure}

Given that the characteristic system of the face $f$ is determined by three parameters, if we consider a symmetry type graph $\mathcal{T(M)}$ described in Figure \ref{classes}, and we remove the edges and semi-edges having colour $k_2$, then we hold a new pregraph $\mathcal{T}_2(\mathcal{M})$ isomorphic to one of those twenty two pregraphs shown in Figure \ref{Im:FacePregraph}.

\begin{figure}[h!]
  \centering
  \includegraphics[scale=0.25]{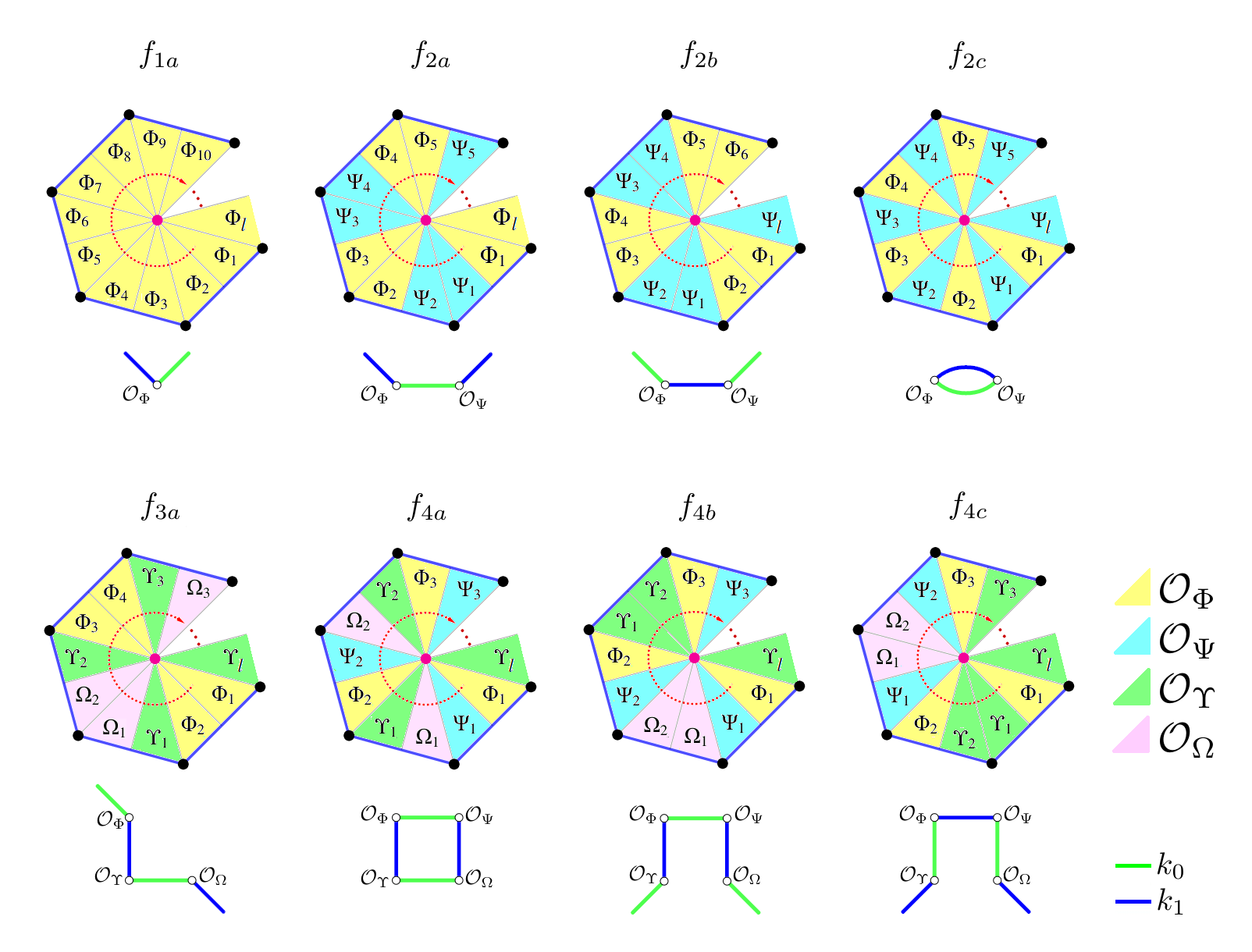}
  \caption{\emph{Sequence follow by the flags in the boundary of a face $f$, from its face type graph of type.}}
  \label{Im:Faces}
\end{figure}

\begin{remark}
The pregraph $\mathcal{T}_2(\mathcal{M})$ is conformed by at most three connected components. Each connected component of $\mathcal{T}_2(\mathcal{M})$, is isomorphic to of one these eight pregraphs shown in Figure \ref{Im:Faces}, which we denote as $f_x$, for some $x$ in the set of index $\{1a, 2a, 2b, 2c, 3a, 4a, 4b, 4c\}$.

\begin{figure}[h!]
  \centering
  \includegraphics[scale=0.25]{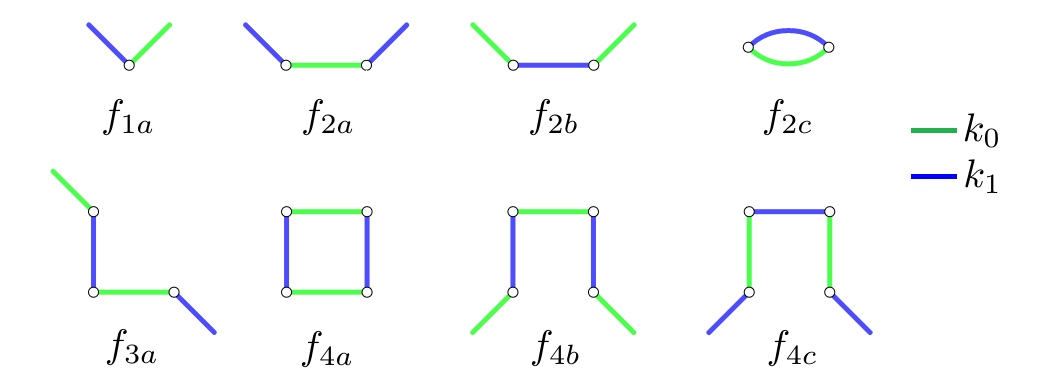}
  \caption{\emph{Pregraphs $f_x$.}}
  \label{Im:FacesPregraph}
\end{figure}
\end{remark}

If we fix a face $f$ of the map $\mathcal{M}$, we hold the cycle $C_{f}$ of $G_{\mathcal{M}}$ such that its vertices are all flags having a vertex in the interior of $f$, and remember that $C_{f}$ is a two colourable alternately cycle, then we can introduce the set of orbits $\mathcal{O}rb(C_{f})$ defined by the action of $Aut(\mathcal{M})$ on $\mathcal{F(M)}$ restricted to the flags that conformed the vertices of $C_{f}$. Using the pregraph definition in Section \ref{sec:pregraph} we hold  that the $\overline{C}_{f}$ is the pregraph of the cycle $C_{f}$ with respect to $\mathcal{O}rb(C_{f})$ induces the following definition.
 
\begin{definition}
Consider the pregraph $\overline{C}_{f}$ which is contained into a connected component of $\mathcal{T}_2\mathcal{(M)}$, by construction $\overline{C}_{f}$ is a pregraph of type $f_x$, for some $x$ in the set of index $\{1a, 2a, 2b,$ $2c, 3a, 4a, 4b, 4c\}$. This connected component is called  \textbf{the face type graph $\mathcal{T}(f)$ of the face $f$}.
\end{definition}

\begin{thm}\label{thm:face_cycle}
Let us considerer a $4$-orbit map and let $f$ be a face of the map. If $f_{3a}$ is the face type graph of $f$ then the boundary of $f$ is conformed by $3n$ edges, for any $n \geq 1$. In other words, $f$ has size $3n$. Moreover, the characteristic system of the face $f$ is $(6n,k_0,k_1)$.
\end{thm}

\begin{proof}
We consider the vertex $i(v)$ and the edge $i(e)$ of the $4$-orbit map $\mathcal{M}$ such that $i(e)$ incidents to $i(v)$ and $i(e)$ belongs to the boundary of the face $f$. Then we denote as $\Phi$ the flag of the map $\mathcal{M}$ conformed by the triplet
\[
\Phi:=(i(v),i(e),f).
\]
Suppose that there are flags $\Upsilon, \Omega$ on the map such that the classes $\mathcal{O}_\Phi$, $\mathcal{O}_\Upsilon$ and $\mathcal{O}_\Omega$ are the vertices of the pregraph $f_{3a}$ (see Figure \ref{face_f3a}-a). We will count the number of elements in the set $\mathcal{O}_{\Phi}^{\langle s_{0},s_{1} \rangle}$ using the face type graph $f_{3a}$.

Let us consider the action of $\langle s_{1},s_{2} \rangle$ on the flags set, by equation (\ref{eq_orbit_face}), the class 
\[
\mathcal{O}_{\Phi}^{\langle s_{0},s_{1} \rangle}=\{\Phi\cdot w : w\in \langle s_{0},s_{1} \rangle\}
\]
contains all the flags having a vertex in the interior of the face $f$, it means
\[
\mathcal{O}_{\Phi}^{\langle s_{0},s_{1} \rangle}=\{\Phi,\Phi^0,\Phi^1, \Phi^{0,1}, \Phi^{1,0}, \Phi^{0,1,0}, \Phi^{1,0,1}, \Phi^{1,0,1,0},... \}
\]

Given that the $\{\mathcal{O}_{\Phi}\}_0$ is a semi-edge of the pregraph $f_{3a}$, then the 0-adjacent flag $\Phi^0$ belongs to  the orbit $\mathcal{O}_{\Phi}$. Then without lost of generality we can rewrite that $\Phi=\Phi_1$ y $\Phi^0=\Phi_2$. Analogously, as the sets $\{\mathcal{O}_{\Phi},\mathcal{O}_{\Upsilon}\}_1$, $\{\mathcal{O}_{\Upsilon},\mathcal{O}_{\Omega}\}_0$, $\{\mathcal{O}_{\Omega}\}_1$, $\{\mathcal{O}_{\Omega},\mathcal{O}_{\Upsilon}\}_0$ and $\{\mathcal{O}_{\Upsilon},\mathcal{O}_{\Phi}\}_1$  are edges of the pregraph $f_{3a}$, then the flags $\Phi^{0,1}$, $\Phi^{0,1,0}$, $\Phi^{0,1,0,1}$, $\Phi^{0,1,0,1,0}$ and $\Phi^{0,1,0,1,0,1}$ are belonged to in the orbit $\mathcal{O}_{\Upsilon}$, $\mathcal{O}_{\Omega}$, $\mathcal{O}_{\Omega}$, $\mathcal{O}_{\Upsilon}$ y $\mathcal{O}_{\Phi}$, respectively. Then we can assume that $\Phi^{0,1}=\Upsilon_1$, $\Phi^{0,1,0}=\Omega_1$, $\Phi^{0,1,0,1}=\Omega_2$, $\Phi^{0,1,0,1,0}=\Upsilon_2$ y $\Phi^{0,1,0,1,0,1}=\Phi_3$. Following with this construction we obtain the finite sequence $\Phi_1,\Phi_2,\Upsilon_1,\Omega_1,\Omega_2,\Upsilon_2,\Phi_3,...,\Phi_{l-1},\Phi_l,\Upsilon_{l-1},\Omega_{l-1},\Omega_l,\Upsilon_l$ (see Figure \ref{face_f3a}-b), where $\Phi_l \in \mathcal{O}_\Phi$, $\Upsilon_l \in \mathcal{O}_\Upsilon$  and $\Omega_l \in \mathcal{O}_\Omega$ for $l \in \{1,2,3,...,n\}$. From this it holds
\[
\begin{array}{ccl}
\mathcal{O}_{\Phi}^{\langle s_{0},s_{1} \rangle}&=&\{\Phi,\Phi^0,\Phi^{0,1}, \Phi^{0,1,0}, \Phi^{0,1,0,1}, \Phi^{0,1,0,1,0}, \Phi^{0,1,0,1,0,1},... \}\\
 & = &\{\Phi_1, \Phi_2, \Upsilon_1, \Omega_1, \Omega_2, \Upsilon_2, \Phi_3,...,\Phi_{l-1},\Phi_l,\Upsilon_{l-1},\Omega_{l-1},\Omega_l,\Upsilon_l\}.         
\end{array}
\]

\begin{figure}[H]
\centering
\begin{tabular}{ccc}
\includegraphics[scale=0.3]{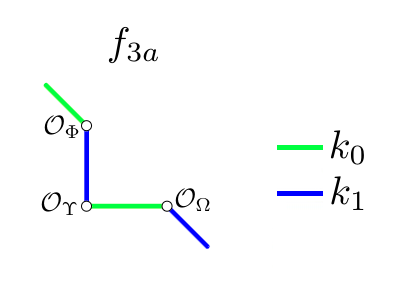}& &\includegraphics[scale=0.3]{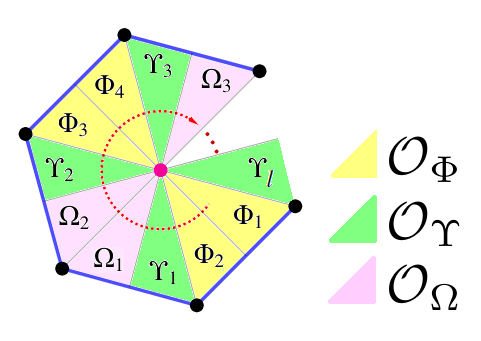}\\
a. Face type graph $f_{3a}$. && b. Orbit $\mathcal{O}_{\Phi}^{\langle s_{0},s_{1} \rangle}$ around \\
   & &  the face $f$.\\
  \end{tabular}
 \caption{\emph{Sequence follow by the flags conforming the boundary of a face $f$, from its face type graph of type $f_{3a}$.}} 
  \label{face_f3a}
\end{figure}

Let us consider that the  face $f$ has $m$ edges on its boundary, we will prove that $m=3n$ for some $n\geq 1$. Let us consider two flags $\Phi_1, \Phi_2 \in \mathcal{O}_{\Phi}^{\langle s_{0},s_{1} \rangle}$, such that they are adjacent by an edge $e$, then we label all edges in the boundary of $f$ in clockwise and with this rewrite  $e=e_1$. By the construction of  $\mathcal{O}_{\Phi}^{\langle s_{0},s_{1} \rangle}$, if $\Delta $ is a flag in $f$ with an edge on $e_l$ such that $1$ is congruent to $l$ modulo $3$, then the flag $\Delta$ is in the class $\mathcal{O}_{\Phi}$, where $l\in\{1,\ldots, m\}$. But if $1$ is not congruent to $l$ modulo $3$, then the flag $\Delta$ is in the class $\mathcal{O}_{\Omega}$ or $\mathcal{O}_{\Upsilon}$. 
%
%$$1 \cong l \quad \text{mod 3}.$$
%%Let us asume that  this is not true, for that, if we call $m$ the size of $f$, let us check $m=3n+r$ with $r$ taking values on $\{1,2\}$.
By division theorem there are two positve integers $n$ and $r$ such that $m=2n+r$ and $r$ taking values on $\{0, 1,2\}$.
If $r=1$, the two flags with edge $e_m$ are belong to the class $\mathcal{O}_{\Phi}$, it means the flags with edges $e_1$ and $e_m$ are in the class $\mathcal{O}_{\Phi}$. However, one  flag with edge $e_m$ must be in the class $\mathcal{O}_{\Omega}$ and the other flag with edge $e_m$ must be on the class $\mathcal{O}_{\Upsilon}$. Thus $r\neq 1$.  Now, if $r=2$, in the edge $e_m$ there is a flag in the class $\mathcal{O}_{\Upsilon}$ and the other in the class $\mathcal{O}_{\Omega}$. Given that flags in the edge $e_{3n+1}$ are in the class $\mathcal{O}_{\Phi}$, then one flag with edge $e_1$ must be in the class $\mathcal{O}_{\Omega}$ and the other flag with edge $e_1$ must be in the class $\mathcal{O}_{\Upsilon}$. However, the both flags with edge $e_1$ are in the class $\mathcal{O}_{\Phi}$. Thus $r\neq 2$ and we conclude that $r=0$. This implies that the number of edges conforming the boundary of $f$ is $3n$, for $n\geq 1$.  Given that there are at least six flags $\Phi_1, \Phi_2, \Upsilon_1, \Omega_1, \Omega_2, \Upsilon_2$ in  $\mathcal{O}_{\Phi}^{\langle s_{0},s_{1} \rangle}$, then it follows that the number of flags in the  class $\mathcal{O}_{\Phi}^{\langle s_{0},s_{1} \rangle}$ is $6n$ for some $n \geq 1$. 
From this, it holds that the characteristic system of the face $f$ is $(6n,k_0,k_1)$.

\end{proof}

If we consider any face $f$ of the $4$-orbit map $\mathcal{M}$ and we suppose that its face type graph $\mathcal{T}(f)$ associated to $f$ is $f_x$, for any $x$ in the set of index $\{1a, 2a, 2b, 2c, 3a, 4a, 4b, 4c\}$, then following the same ideas in the proof of Theorem \ref{thm:face_cycle}, it is easy to find the number of edges conforming the boundary of $f$. These results are collected in Table \ref{Table:FaceType}. In Figure \ref{Im:Faces} are represented the face type graphs.

\begin{table}[H]
\begin{center}
\begin{tabular}{|c|c|c|c|}
\hline 
 Face type graph & Size of $f$& & Characteristic System\\
 of the face $f$ & & & of the face $f$\\
\hline 
$f_{1a}$ & $n$ & $n \geq 3$& $(2n,k_0,k_1)$ \\
\hline 
$f_{2a}$ & $2n$ & $n \geq 2$& $(4n,k_0,k_1)$ \\
\hline 
$f_{2b}$ & $2n$ & $n \geq 2$ & $(4n,k_0,k_1)$ \\
\hline 
$f_{2c}$ & $n$ & $n \geq 3$& $(2n,k_0,k_1)$ \\
\hline 
$f_{3a}$ & $3n$ & $n \geq 1$& $(6n,k_0,k_1)$  \\
\hline 
$f_{4a}$ & $2n$ & $n \geq 2$& $(4n,k_0,k_1)$  \\
\hline 
$f_{4b}$ & $4n$ & $n \geq 1$& $(8n,k_0,k_1)$ \\
\hline 
$f_{4c}$ & $4n$ & $n \geq 1$&$(8n,k_0,k_1)$  \\
\hline 
\end{tabular}
\caption{\emph{Size and characteristic system of a face  $f$ from its face type graph.}} 
\label{Table:FaceType}
\end{center}
\end{table}

\begin{remark}\label{Re:FaceType}
Let $f_1, f_2$ be faces of a $4$-orbit map $\mathcal{M}$ and let $C_{f_1},C_{f_2}$ be the cycles of the flag graph $G_{\mathcal{M}}$, associated to the faces $f_1$ and $f_2$, respectively. Suppose that the pregraphs $\overline{C}_{f_1}$ and $\overline{C}_{f_2}$ are contained into some connected component of $\mathcal{T}_{2}(\mathcal{M})$. Given that the automorphism group $Aut(\mathcal{M})$ acts freely on the set of flags $\mathcal{F(M)}$, then the faces $f_1$ and $f_2$ have the same number of edges in its boundary, and their characteristic systems is
\[
(2m_{f_1},k_0,k_1)=(2m_{f_2},k_0,k_1).
\]
However, if the pregraphs $\overline{C}_{f_1}$ and $\overline{C}_{f_2}$ belong to different connected component of $\mathcal{T}_2(\mathcal{M})$ but they are isomorphic, then the size of the faces $f_1$ and $f_2$ are multiples of $s$, for any $s\in\mathbb{N}$.
If $l$ is number of connected component of $\mathcal{T}_2(\mathcal{M})$, for any $l\in\{1,2,3\}$, then there are $l$ values to the size of the faces of $\mathcal{M}$. Hence, there are $n_i$ positive integers,  with $i\in\{1,\ldots,l\}$ such that the characteristic system of any face of the map is either $(2n_1,k_0,k_1),\ldots,(2n_l,k_0,k_1)$.

\end{remark}

\subsection{Main consequence}

With the elements introduced until this point, we shall study the $4$-orbit maps having symmetry type graph  $4_A$, and we shall summarize through a table the feasible values taken by the degree of the vertices and the appropriate number of edges in the boundary of each face of the $4$-orbit map.

\begin{thm}\label{thm_General}

If the symmetry type graph of a $4$-orbit map $\mathcal{M}$ is $4_{A}$, then 
\begin{enumerate}
\item The pregraph $\mathcal{T}_0(\mathcal{M})$ has only one  connected component isomorphic to $v_{4b}$. If $i(v)$ is a vertex of $\mathcal{M}$, then there is a positive integer $n$ such that the degree of $i(v)$ is $4n$, the characteristic system of $i(v)$ is $(8n,k_1,k_2)$, and its vertex type graph is $v_{4b}$.

\item The pregraph $\mathcal{T}_2(\mathcal{M})$ is conformed by two  connected components isomorphic to $f_{2b}$. If $f$ is a face of $\mathcal{M}$, then there are positive integers $m,n$ such that the number of edges in the boundary of $f$ is either $2n$ or $2m$, the characteristic system of $f$ is either $(4n,k_0,k_1)$ or $(4m,k_0,k_1)$, and its face type pregraph is $f_{2b}$.
\end{enumerate}
\end{thm}

\begin{proof}
If we remove the edges and semi-edges of $4_A$ having colour $k_0$, then the new pregraph $\mathcal{T}_0(\mathcal{M})$ is conformed by  a connected component isomorphic to $v_{4b}$ (see Figures \ref{Im:VertexPregraph} and \ref{Im:VerticesPregraph}). This implies that for each vertex $i(v)$ of the map $\mathcal{M}$ it has vertex type graph $v_{4b}$ and characteristic system $(8n,k_1,k)$ (see Table \ref{Table:VertexType}). This properties are summarized in Table \ref{Table:Type4A}.

Analogously, if we remove the edges and semi-edges of $4_A$ having colour $k_2$, then the new pregraph $\mathcal{T}_2(\mathcal{M})$ is conformed by  two connected components isomorphic to $f_{2a}$ (see Figures \ref{Im:FacePregraph} and \ref{Im:FacesPregraph}). This implies that for each face $f$ of the map $\mathcal{M}$ it has face type graph $f_{2a}$ and there are positive integers $m_1,m_2$ such that its characteristic system is either  $(4m_1,k_0,k,1)$ or $(4m_2,k_0,k,1)$ (see Table \ref{Table:FaceType}). This properties are summarized in Table \ref{Table:Type4A}.
\end{proof}

\begin{table}[h]
\begin{center}
\begin{tabular}{|c|c|c|c|c|c|}
\hline 

 Pregraph &  Number of & Vertex type & \multicolumn{2}{c|} {Degree of the}& Characteristic\\

$\mathcal{T}_0(\mathcal{M})$&  connected component & graph of the &\multicolumn{2}{c|}{ vertex $i(v)$}& System of\\

&  of $\mathcal{T}_0(\mathcal{M})$&{\small vertex $i(v)$} & \multicolumn{2}{c|}{}& the vertex $i(v)$\\ 
\hline 
$4_{(A)_{0}}$ & $1$ & $v_{4b}$ & $4n$  & $n \geq 1$& $(8n,k_1,k_2)$\tabularnewline
\hline
 Pregraph &  Number of &  Face type & \multicolumn{2}{c|}{ Size of the} & Characteristic\\

 $\mathcal{T}_2(\mathcal{M})$ & connected component &{\small graph of the} &\multicolumn{2}{c|}{ face $f$}& System of\\

& of $\mathcal{T}_2(\mathcal{M})$& face $f$ & \multicolumn{2}{c|}{}& the face $f$\\
\hline 
\multirow{2}{*}{$4_{(A)_{2}}$} & \multirow{2}{*}{2} & $f_{2b}$ & $2m_1$ & $m_1 \geq 2$& $(4m_1,k_0,k_1)$\tabularnewline
%\cline{3-5} 
 &  & $f_{2b}$ & $2m_2$ &  $m_2 \geq 2$& $(4m_2,k_0,k_1)$\tabularnewline
\hline 
\end{tabular}
\caption{\emph{Properties for $4$-orbit maps with symmetry type graph $4_A$.}} 
  \label{Table:Type4A}
\end{center}
\end{table}

Following the same ideas that in the proof of the Theorem \ref{thm_General} for any other of the twenty one possibles symmetry type graphs associated to the $4$-orbit maps, the characterization in terms of number of connected components of $\mathcal{T}_{0}(\mathcal{M})$, vertex type graph, degree of a vertex and characteristic system of a vertex are given in the Table \ref{Table:Pregrafos_T0}. Respectively, the characterization in terms of number of connected components of $\mathcal{T}_{2}(\mathcal{M})$, face type graph, size of a face and characteristic system of a face are given in the Table \ref{Table:Pregrafos_T2}.
From the definition of dual map and the Proposition \ref{prop:Hubard13} it follows that:

\begin{corolario}
If $\mathcal{M}$ is a $4$-orbit map with symmetry type graph $\mathcal{T(M)}$  then  

\begin{enumerate}
\item The pregraph$\mathcal{T}_0(\mathcal{M})$  is isomorphic to the pregraph $\mathcal{T}_2(\mathcal{M^*})$.

\item The pregraph  $\mathcal{T}_2(\mathcal{M})$ is isomorphic to the pregraph $\mathcal{T}_0(\mathcal{M^*})$.
\end{enumerate}

\end{corolario}

\newpage

\begin{table}[h!]
\begin{center}
\begin{tabular}{|c|c|c|c|c|c|}
\hline 

 Pregraph &  Number of & Vertex type & \multicolumn{2}{c|}{ Degree of the} &  Characteristic\\

 $\mathcal{T}_0(\mathcal{M})$& connected component & graph of the &\multicolumn{2}{c|}{ vertex $i(v)$}& System of\\

& of $\mathcal{T}_0(\mathcal{M})$& vertex $i(v)$ & \multicolumn{2}{c|}{}& the vertex $i(v)$\\

\hline 
$4_{(A)_{0}}$ & $1$ & $v_{4b}$ & $4n$  & $n \geq 1$& $(8n,k_1,k_2)$\\
\hline
\multirow{2}{*}{$4_{(Ad)_{0}}$} & \multirow{2}{*}{2} & $v_{2b}$ & $2n_1$  &  $n_1 \geq 2$ & $(4n_1,k_1,k_2)$\tabularnewline
%\cline{3-5} 
 &  & $v_{2b}$ & $2n_2$  &  $n_2\geq 2$  & $(4n_2,k_1,k_2)$\tabularnewline

\hline 
$4_{(Ap)_{0}}$ & 1 &$v_{4b}$ & $4n$  &  $n \geq 1$& $(8n,k_1,k_2)$\tabularnewline
\hline 

$4_{(B)_{0}}$ & 1 & $v_{4c}$ & $4n$  &  $n \geq 1$& $(8n,k_1,k_2)$\tabularnewline
\hline 

\multirow{3}{*}{$4_{(Bd)_{0}}$} & \multirow{3}{*}{3} & $v_{2b}$ & $2n_1$  &  $n_1 \geq 2$& $(4n_1,k_1,k_2)$\tabularnewline
%\cline{3-5} 
 &  & $v_{1a}$ & $n_2$  &  $n_2 \geq 3$& $(2n_2,k_1,k_2)$\tabularnewline
%\cline{3-5} 
 &  & $v_{1a}$ & $n_3$  &  $n_3 \geq 3$& $(2n_3,k_1,k_2)$\tabularnewline
\hline 

$4_{(Bp)_{0}}$ & 1 & $v_{4c}$ & $4n$  &  $n \geq 1$& $(8n,k_1,k_2)$\tabularnewline
\hline 

$4_{(C)_{0}}$ & 1 & $v_{4a}$ & $2n$  &  $n \geq 2$& $(4n,k_1,k_2)$\tabularnewline
\hline 

\multirow{2}{*}{$4_{(Cd)_{0}}$} & \multirow{2}{*}{2} & $v_{2b}$ & $2n_1$  &  $n_1 \geq 2$& $(4n_1,k_1,k_2)$\tabularnewline
%\cline{3-5} 
 &  & $v_{2b}$ & $2n_2$  &   $n_2\geq 2$& $(4n_2,k_1,k_2)$\tabularnewline
\hline 

$4_{(Cp)_{0}}$ & 1 & $v_{4a}$ & $2n$  &  $n \geq 2$& $(4n,k_1,k_2)$\tabularnewline
\hline 

\multirow{2}{*}{$4_{(D)_{0}}$} & \multirow{2}{*}{2} & $v_{3a}$ & $3n_1$  &  $n_1 \geq 1$& $(6n_1,k_1,k_2)$\tabularnewline
%\cline{3-5} 
 &  & $v_{1a}$ & $n_2$  &  $n_2\geq 3$& $(2n_2,k_1,k_2)$\tabularnewline
\hline 

$4_{(Dd)_{0}}$ & 1 & $v_{4c}$ & $4n$  &  $n \geq 1$& $(8n,k_1,k_2)$\tabularnewline
\hline 

\multirow{2}{*}{$4_{(Dp)_{0}}$} & \multirow{2}{*}{2} & $v_{3a}$ & $3n_1$  &  $n_1 \geq 1$& $(6n_1,k_1,k_2)$\tabularnewline
%\cline{3-5} 
 &  & $v_{1a}$ & $n_2$  &  $n_2 \geq 3$& $(2n_2,k_k,k_2)$\tabularnewline
\hline 

$4_{(E)_{0}}$ & 1 & $v_{4b}$ & $4n$  &  $n \geq 1$& $(8n,k_1,k_2)$\tabularnewline
\hline 

$4_{(Ed)_{0}}$ & 1 & $v_{4a}$ & $2n$  &  $n \geq 2$& $(4n,k_1,k_2)$\tabularnewline
\hline 

$4_{(Ep)_{0}}$ & 1 & $v_{4b}$ & $4n$  &  $n \geq 1$& $(8n,k_1,k_2)$\tabularnewline
\hline 

\multirow{2}{*}{$4_{(F)_{0}}$} & \multirow{2}{*}{2} & $v_{2a}$ & $2n_1$  &  $n_1 \geq 2$& $(4n_1,k_1,k_2)$\tabularnewline
%\cline{3-5} 
 &  & $v_{2a}$ & $2n_2$  &  $n_2 \geq 2$& $(4n_2,k_1,k_2)$\tabularnewline
\hline 

$4_{(G)_{0}}$ & 1 & $v_{4a}$ & $2n$  &  $n \geq 2$& $(4n,k_1,k_2)$\tabularnewline
\hline
 
\multirow{2}{*}{$4_{(Gd)_{0}}$} & \multirow{2}{*}{2} & $v_{2c}$ & $n_1$  &  $n_1 \geq 3$& $(2n_1,k_1,k_2)$\tabularnewline
%\cline{3-5} 
 &  & $v_{2c}$ & $n_2$  &  $n_2 \geq 3$& $(2n_2,k_1,k_2)$\tabularnewline
\hline 

$4_{(Gp)_{0}}$ & 1 & $v_{4a}$ & $2n$  &  $n \geq 2$& $(4n,k_1,k_2)$\tabularnewline
\hline 

$4_{(H)_{0}}$ & 1 & $v_{4c}$ & $4n$  &  $n \geq 1$& $(8n,k_1,k_2)$\tabularnewline
\hline 

\multirow{2}{*}{$4_{(Hd)_{0}}$} & \multirow{2}{*}{2} & $v_{2c}$ & $n_1$  &  $n_1 \geq 3$& $(2n_1,k_1,k_2)$\tabularnewline
%\cline{3-5} 
 &  & $v_{2a}$ & $2n_2$  &  $n_2\geq 2$& $(4n_2,k_1,k_2)$\tabularnewline
\hline 

$4_{(Hp)_{0}}$ & 1 & $v_{4c}$ & $4n$  &  $n \geq 1$& $(8n,k_1,k_2)$\tabularnewline
\hline 
\end{tabular}
 \caption{\emph{Properties for $4$-orbit maps from its pregraph $\mathcal{T}_{0}(\mathcal{M})$.}} 
\label{Table:Pregrafos_T0}
\end{center}
\end{table}

\begin{table}[H]
\begin{center}
\begin{tabular}{|c|c|c|c|c|c|}
\hline 
Pregraph &  Number of &  Face type & \multicolumn{2}{c|}{Size of the} &  Characteristic\\

$\mathcal{T}_2(\mathcal{M})$& connected component & graph of the &\multicolumn{2}{c|}{ face $f$}& System of\\

& of  $\mathcal{T}_2(\mathcal{M})$ & face $f$ & \multicolumn{2}{c|}{}& the face $f$\\
\hline 
\multirow{2}{*}{$4_{(A)_{2}}$} & \multirow{2}{*}{2} & $f_{2b}$ & $2m_1$ & $m_1 \geq 2$& $(4m_1,k_0,k_1)$\tabularnewline
%\cline{3-5} 
 &  & $f_{2b}$ & $2m_2$ &  $m_2 \geq 2$& $(4m_2,k_0,k_1)$\tabularnewline
\hline 

$4_{(Ad)_{2}}$ & 1 & $f_{4b}$ & $4m$ &  $m \geq 1$& $(8m,k_0,k_1)$\tabularnewline
\hline
 
$4_{(Ap)_{2}}$ & 1 & $f_{4b}$ & $4m$ &  $m \geq 1$& $(8m,k_0,k_1)$\tabularnewline
\hline 

\multirow{3}{*}{$4_{(B)_{2}}$} & \multirow{3}{*}{3} & $f_{2b}$ & $2m_1$ &  $m_1 \geq 2$& $(4m_1,k_0,k_1)$\tabularnewline
%\cline{3-5} 
 &  & $f_{1a}$ & $m_2$ &  $m_2 \geq 3$& $(2m_2,k_0,k_1)$\tabularnewline
%\cline{3-5} 
 &  & $f_{1a}$ & $m_3$ &  $m_3 \geq 3$& $(2m_3,k_0,k_1)$\tabularnewline
\hline
 
$4_{(Bd)_{2}}$ & 1 & $f_{4c}$ & $4m$ &  $m \geq 1$& $(8m,k_0,k_1)$\tabularnewline
\hline
 
$4_{(Bp)_{2}}$ & 1 & $f_{4c}$ & $4m$ &  $m \geq 1$& $(8m,k_0,k_1)$\tabularnewline
\hline 

\multirow{2}{*}{$4_{(C)_{2}}$} & \multirow{2}{*}{2} & $f_{2b}$ & $2m_1$ &  $m_1 \geq 2$& $(4m_1,k_0,k_1)$\tabularnewline
%\cline{3-5} 
 &  & $f_{2b}$ & $2m_2$ &  $m_2 \geq 2$& $(4m_2,k_0,k_1)$\tabularnewline
\hline 

$4_{(Cd)_{2}}$ & 1 & $f_{4a}$ & $2m$ &  $m \geq 2$& $(4m,k_0,k_1)$\tabularnewline
\hline

$4_{(Cp)_{2}}$ & 1 & $f_{4a}$ & $2m$ &  $m \geq 2$& $(4m,k_0,k_1)$\tabularnewline
\hline 

\multicolumn{1}{|c|}{$4_{(D)_{2}}$} & \multicolumn{1}{c|}{1} & $f_{4c}$ & $4m$ &  $m \geq 1$& $(8m,k_0,k_1)$\tabularnewline
\hline

\multirow{2}{*}{$4_{(Dd)_{2}}$} & \multirow{2}{*}{2} & $f_{3a}$ & $3m_1$ &  $m_1 \geq 1$& $(6m_1,k_0,k_1)$\tabularnewline
%\cline{3-5} 
 &  & $f_{1a}$ & $m_2$ &  $m_2 \geq 3$& $(2m_2,k_0,k_1)$\tabularnewline
\hline 

\multirow{2}{*}{$4_{(Dp)_{2}}$} & \multirow{2}{*}{2} & $f_{3a}$ & $3m_1$ &  $m_1 \geq 1$& $(6m_1,k_0,k_1)$\tabularnewline
%\cline{3-5} 
 &  & $f_{1a}$ & $m_2$ &  $m_2 \geq 3$& $(2m_2,k_0,k_1)$\tabularnewline
\hline

$4_{(E)_{2}}$ & 1 & $f_{4a}$ & $2m$ &  $m \geq 2$& $(4m,k_0,k_1)$\tabularnewline
\hline

$4_{(Ed)_{2}}$ & 1 & $f_{4b}$ & $4m$ &  $m \geq 1$& $(8m,k_0,k_1)$\tabularnewline
\hline
 
$4_{(Ep)_{2}}$ & 1 & $f_{4b}$ & $4m$ &  $m \geq 1$& $(8m,k_0,k_1)$\tabularnewline
\hline
 
\multirow{2}{*}{$4_{(F)_{2}}$} & \multirow{2}{*}{2} & $f_{2a}$ & $2m_1$ &  $m_1 \geq 2$& $(4m_1,k_0,k_1)$\tabularnewline
%\cline{3-5} 
 &  & $f_{2a}$ & $2m_2$ &  $m_2 \geq 2$& $(4m_2,k_0,k_1)$\tabularnewline

\hline 
\multirow{2}{*}{$4_{(G)_{2}}$} & \multirow{2}{*}{2} & $f_{2c}$ & $m_1$ &  $m_1 \geq 3$& $(2m_1,k_0,k_1)$\tabularnewline
%\cline{3-5} 
 &  & $f_{2c}$ & $m_2$ &  $m_2 \geq 3$& $(2m_2,k_0,k_1)$\tabularnewline
\hline
 
$4_{(Gd)_{2}}$ & 1 & $f_{4a}$ & $2m$ &  $m \geq 2$& $(4m,k_0,k_1)$\tabularnewline
\hline 

$4_{(Gp)_{2}}$ & 1 & $f_{4a}$ & $2m$ &  $m \geq 2$& $(4m,k_0,k_1)$\tabularnewline
\hline 

\multirow{2}{*}{$4_{(H)_{2}}$} & \multirow{2}{*}{2} & $f_{2c}$ & $m_1$ &  $m_1 \geq 3$& $(2m_1,k_0,k_1)$\tabularnewline
%\cline{3-5} 
 &  & $f_{2a}$ & $2m_2$ &  $m_2 \geq 2$& $(4m_2,k_0,k_1)$\tabularnewline
\hline

$4_{(Hd)_{2}}$ & 1 & $f_{4c}$ & $4m$ &  $m \geq 1$& $(8m,k_0,k_1)$\tabularnewline
\hline
 
$4_{(Hp)_{2}}$ & 1 & $f_{4c}$ & $4m$ &  $m \geq 1$& $(8m,k_0,k_1)$\tabularnewline
\hline 
\end{tabular}
\caption{\emph{Properties for $4$-orbit maps from its pregraph $\mathcal{T}_{2}(\mathcal{M})$.}} 
\label{Table:Pregrafos_T2}
\end{center}
\end{table}

%%%%%%%%%%%%%%%%%%%%%%%%%%%%%%%%%%%%%%%%%%%
%%%%%%%%%%%%%%%%%%%%%%%%%%%%%%%%%%%%%%%%%%%%%%%

\section*{Acknowledgements}
The second author was partially supported by UNIVERSIDAD NACIONAL DE COLOMBIA, SEDE MANIZALES. Camilo Ram\'irez Maluendas has dedicated this work to his beautiful family: Marbella and Emilio, in appreciation of their love and support.

The third author wishes to thank her colleagues Alexander and Camilo for their valuable teachings, dedication and support. In addition, the author thanks to the Fundaci\'on Universitaria Konrad Lorenz for the support.   
%%%%%%%%%%%%%%%%%%%%%%%%%%%
%%%%%%%%%%%%%%%%%%%%%%%%%%%%%%%
%%Ac\'a comienza lo chido

%%%%%%%%%%%%%%%%%%%%%%%%%%%%%%%%%%%%%%%%%%%%%%%%%%%%%%%%%%%%%%%%%%%%%%%%%%%%%%%%%%%%%%%%%%%%%%%%%%%%%%%%%%%%%%%%%%%%%%%%%%%%%%%%%%%%%%%%%%%%%%%%%%%%%%%%%%%%%%%%%%

\end{document}